\documentclass[DIV=13]{scrartcl}
\usepackage[sb]{libertine} %
\usepackage[T1]{fontenc}
\usepackage{textcomp}
\usepackage[varqu,varl]{zi4}%
\usepackage[amsthm,upint]{libertinust1math} %
\usepackage[scr=boondoxo,bb=boondox]{mathalpha} %
\usepackage{bm}
\usepackage{mathtools}
\usepackage{todonotes}
\usepackage{booktabs}
\usepackage{enumitem}
\usepackage{placeins}
\usepackage{colortbl}
\usepackage{tikz}
\usepackage{bm}
\usepackage{biblatex}
\AtEveryBibitem{%
  \clearfield{urlyear}%
  \clearfield{urldate}%
}%
\AtEveryCitekey{%
  \clearfield{urlyear}%
  \clearfield{urldate}%
}%
\usepackage{xcolor}
\usepackage{float}
\RequirePackage{siunitx}
\RequirePackage[algoruled,norelsize, linesnumbered, lined, commentsnumbered]{algorithm2e}
\newcommand{\blu}[1]{{\color{myblue} \noindent #1}}

\usepackage{mycommands}
\usepackage[format=plain,labelfont=bf]{caption}
\usepackage[pdftex,colorlinks=true,linkcolor={blue!50!black},citecolor={red!50!black},urlcolor=blue]{hyperref}

\renewcommand{\vec}[1]{\mathbf{#1}}
\newcommand{\fvec}[1]{\bm{#1}}

\definecolor{myblue}{RGB}{0 83 139}
\definecolor{myred}{RGB}{114 16 69}
\definecolor{mygreen}{RGB}{0 94 0}
\newcommand{\bestf}[1]{{\color{myblue}\bfseries \noindent #1}}
\newcommand{\bests}[1]{{\color{mygreen}\bfseries \noindent #1}}
\newcommand{\bestt}[1]{{\color{myred}\bfseries \noindent #1}}
\newcommand\restrict[1]{\raisebox{-.5ex}{$|$}_{#1}}

\begin{document}
\date{}
\title{DNN-MG: A Hybrid Neural Network/Finite Element Method with Applications to 3D Simulations of the Navier-Stokes Equations}
\author{Nils Margenberg\thanks{Corresponding author}
  \thanks{Helmut Schmidt University,
    Faculty of Mechanical and Civil Engineering,
    Holstenhofweg 85,
    22043 Hamburg,
    Germany,
    \href{mailto:margenbn@hsu-hh.de}{\texttt{margenbn@hsu-hh.de}}}
  \and Robert Jendersie\thanks{
    University of Magdeburg,
    Institute for Analysis and Numerics,
    Universit\"atsplatz 2,
    39104 Magdeburg,
    Germany,
    \href{mailto:robert.jendersie@ovgu.de}{\texttt{robert.jendersie@ovgu.de}}}
  \and Christian Lessig\thanks{
    University of Magdeburg,
    Institute for Simulation and Graphics,
    Universit\"atsplatz 2,
    39104 Magdeburg,
    Germany,
    \href{mailto:christian.lessig@ovgu.de}{\texttt{christian.lessig@ovgu.de}}}
  \and Thomas Richter\thanks{
    University of Magdeburg,
    Institute for Analysis and Numerics,
    Universit\"atsplatz 2,
    39104 Magdeburg,
    Germany,
    \href{mailto:thomas.richter@ovgu.de}{\texttt{thomas.richter@ovgu.de}}}
}

\maketitle
\vspace{-7ex}
\begin{abstract}
We extend and analyze the deep neural network multigrid solver (DNN-MG) for the Navier-Stokes equations in three dimensions. The idea of the method is to augment a finite element simulation on coarse grids with fine scale information obtained using deep neural networks.
The neural network operates locally on small patches of grid elements. The
local approach proves to be highly efficient, since the network can be kept (relatively) small and since it can be applied in parallel on all grid patches. However, the main advantage of the local approach is the inherent  generalizability of the method. Since the network only processes data of small sub-areas, it never ``sees'' the global problem and thus does not learn false biases.
We describe the method with a focus on the interplay between the finite element method and deep neural networks. Further, we demonstrate with numerical examples the excellent efficiency of the hybrid approach, which allows us to achieve very high accuracy with a coarse grid and thus reduce the computation time by orders of magnitude.
\end{abstract}

\section{Introduction}
\label{sec:intro}
Accurate flow simulations remain a challenging task. The success of deep neural
networks in machine translation, computer vision, and many other fields has lead to a growing (and renewed) interest to apply neural
networks to problems in computational science and engineering, including flow simulations. The combination
of classical finite element approximation techniques with deep neural networks
adds new aspects to pure numerics-oriented approaches. 

For fluid flows, an accurate simulation, especially in 3D and at higher Reynolds numbers, is still challenging and classical methods reach their limits when very high accuracy is required but also fast results. Although the finite element method (FEM) is
highly efficient and established for the discretization of the Navier-Stokes
equations, fundamental problems, such as the resolution of fine structures or a
correct information transport between scales, are still not
fully solved. 
An accurate approximation of the Navier-Stokes equations is, furthermore, hampered by additional local features under mesh refinement that appear for the nonlinear problem so that also linear solvers with optimal complexity are no remedy.
In this work, we investigate first steps for how a hybrid neural network enriched FEM simulation for 3D flow that can improve the computational speed and provide accurate approximations.

(Adaptive) finite elements can be considered as a hierarchical method, where a hierarchy of finite element meshes and spaces can be constructed (adaptively) to yield approximations of increasing accuracy. This hierarchical setup is then used within multigrid methods for the efficient solution of the algebraic systems. Here, we connect such a hierarchical finite element approach with neural networks in a hybrid setting: while coarse mesh levels are handled in the traditional finite element multigrid way, updates one finer mesh levels are learned by the neural network. They are hence used whenever a full resolution of the effects does not seem possible or efficient.

\begin{figure}
  \centering
  \includegraphics[width=\linewidth]{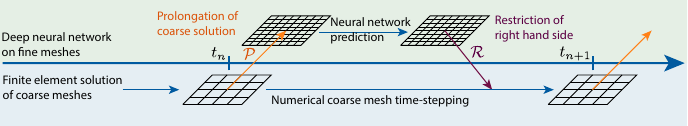}
  \caption{
    The idea of DNN-MG is to use a fast coarse grid solver (e.\,g.\ a multigrid) to get a quick approximation of the solution. This obtained solution
    is then prolongated to a finer level $L+J$, where a neural network predicts
    corrections to the solution. These corrections are incorporated into the
    time evolution through the time stepping and enter through the right-hand
    side. }\label{fig:g-abstract}
\end{figure}

We call this approach the Deep Neural Network Multigrid Solver (DNN-MG) as it is based on hierarchies of meshes and functions spaces and combines the tools of multigrid solvers with deep neural networks (cf.\ Figure~\ref{fig:g-abstract}). However, the neural networks are used to predict updates to the nonlinear problem and hence, the approach can be used for an upsampling of any kind of finite element solution (and also finite volume or finite difference)  towards a representation on a finer mesh. Even though we have tightly embedded the method in a geometric multigrid method, the coarse grid problem can be approximated with any solver, e.g., a Newton-Krylov space method, or pressure-projection type flow solvers. Also, the prolongation onto the fine grid on which the neural network acts need not be a classical multigrid transfer. Indeed, local enrichment is also possible without a global fine grid. Instead the solution could, for example, be corrected in a higher order space.

In this work, we demonstrate the efficiency, generalizability, and scalability of DNN-MG for 3D simulations. We first train the neural network on data from
the classical channel flow with one circular obstacle. To analyze the generalization
capability of DNN-MG, we test the network on channel flows with one or two
obstacles at different Reynolds numbers. The obstacles have an elliptical
cross-section with varying eccentricities. The obtained solutions as well as the
lift and drag functionals demonstrate that DNN-MG obtains considerably better
accuracy than a coarse  solution while taking less
than $3\%$ of the additional computation time required by a full solution on a twice refined mesh.
Therefore, DNN-MG offers a speedup by a factor of \(35\). DNN-MG's efficiency is
evident, as it requires only double the time of a coarse mesh solution on level
$L$, yet offers substantial improvements in the overall solution quality.

The paper is organized as follows. In the Section~\ref{sec:related-work}, we review related work
on using neural networks for the simulation of partial differential equations.
Sections~\ref{sec:num-meth} provides a recap of the solution of the Navier-Stokes equations using the
geometric multigrid method, while Section~\ref{sec:nn} introduces the neural networks used
later in our numerical experiments. In Section~\ref{sec:dnnmg}, we present the 3D deep neural
network multigrid solver, discussing its design in a general form, which makes
it applicable to other problems.
Finally, we present our numerical results in Section~\ref{sec:num-example}.

\section{Related works}
\label{sec:related-work}
We discuss different approaches to the simulation of physical systems using deep
neural networks.

The most direct approach is to represent the solution of the PDE by a (deep) neural network. In~\cite{lagarisArtificialNeuralNetworks1998} already, a least squares formulation has been employed to minimize the residual in the training process. The idea has recently been used again in the Deep Ritz Method~\cite{eDeepRitzMethod2018} and a variety of similar approaches have emerged in the last years and coined Physics Informed Neural Network (PINNs) in the literature~\cite{luDeepXDEDeepLearning2021}. For a comprehensive review of PINN and related approaches, we refer to~\cite{cuomoScientificMachineLearning2022}. Such approaches use an existing partial differential equation and can be considered data free since no explicit training data is required and instead the equations itself is incorporated into the loss function. They are, however, static in that the solution to one specific problem instance is learned. PINNs show superiority for instance in very high dimensional problems~\cite{luPrioriGeneralizationAnalysis2021} but they do not reach the efficiency of established methods when standard problems, e.\,g.\ from elasticity or fluid dynamics, are considered (cf.~\cite{grossmannCanPhysicsInformedNeural2023} for case studies on the Poisson, Allen-Cahn and Schrödinger equation). The idea of PINNs has also been extended to learning solution operators~\cite{luLearningNonlinearOperators2021,chenUniversalApproximationNonlinear1995,kutyniok2022atheoretical}. The main advantage of this approach is that once the solution operator is trained, it can be applied to other settings. This enables the application to inverse problems~\cite{molinaroNeuralInverseOperators2023,kaltenbachSemisupervisedInvertibleNeural2023,caoResidualbasedErrorCorrection2023} or optimal control~\cite{margenbergOptimalDirichletBoundary2023}. Inspired by classical time stepping schemes, another approach to a neural network-based simulation of (evolutionary) partial differential equations is to consider these as time series and predict the next state of the system based on the previous ones. With this approach, network architecture for sequential data can be used, in particular recurrent neural networks and their extensions such as Long-Short-Term-Memory units (LSTM)~\cite{hochreiterLongShorttermMemory1997} and Gated Recurrent Units (GRUs)~\cite{choLearningPhraseRepresentations2014} that we also used for the DNN-MG Navier-Stokes solver in 2D~\cite{margenbergStructurePreservationDeep2021,margenbergNeuralNetworkMultigrid2022}.

At the other end of the spectrum are simulation techniques that are purely data-based. For example, in~\cite{eichingerStationaryFlowPredictions2021,eichinger2022} flow fields are learned using convolutional neural networks based on the flow geometry and using techniques from image processing. Recent work~\cite{pathakFourCastNetGlobalDatadriven2022,biPanguWeather3DHighResolution2022,lamGraphCastLearningSkillful2022} uses purely data-driven neural network models trained on historical weather measurements and these perform on par with the most sophisticated weather forecasting models based on partial differential equations. The authors of~\cite{genevaTransformersModelingPhysical2022} use transformer models for the prediction of dynamical systems. In the context of dynamical systems, the combination of data assimilation, uncertainty quantification, and machine learning techniques has gained significant interest, we refer to~\cite{chengMachineLearningData2023} for an in-depth review. For a recent overview on different neural network based approaches to approximate partial differential equations and related inverse problems we refer to~\cite{Tanyu2023}.

One of the central open questions in the current literature is
to what extent and how physical constraints or existing knowledge should be used for the training of a neural
network, see also~\cite{ghattasLearningPhysicsbasedModels2021}. The Deep Neural Network Multigrid Solver (DNN-MG) is a hybrid approach that combines the model based finite element method with deep neural networks trained on example data that includes model knowledge through the fine mesh residual. Although the hybrid DNN-MG approach integrates model knowledge, DNN-MG is therefore not a PINN. The general setup of DNN-MG aims for generalization, and it is closer to learning a numerical solver than to learning a specific solution.

The close incorporation of the deep neural network into the finite element method, e.\,g.\ that we learn fine-scale solution coefficients, gives rise to analytical tools that can be used in first steps of error estimation~\cite{minakowskiPrioriPosterioriError2023}. The techniques are used in the error estimation for the DNN-MG method~\cite{kapustsinHybridFiniteElement,kapustsinErrorAnalysisHybrid2023}.

In our first work on DNN-MG~\cite{margenbergNeuralNetworkMultigrid2022} we introduced a method which
combines a coarse mesh solution with ANNs but only considered 2D
simulations. We also adapted the method in 2D to ensure divergence-freedom by
construction through a network that predicts the streamfunction~\cite{margenbergStructurePreservationDeep2021}. In later studies we showed that
the corrections added by the neural network improve the guesses in the Newton
method, which leads to a reduction in wall-time compared to the coarse grid
solutions~\cite{margenbergDeepNeuralNetworks2021}. Similar to the ideas of
DNN-MG, the authors of~\cite{delaraAcceleratingHighOrder2022,manriquedelaraAcceleratingHighOrder2023} use numerical
solutions of higher order to train a neural network which generates corrective
forcing terms.
Similarly, in~\cite{fabraFiniteElementApproximation2022} low-order simulations
of wave equations
are enhanced using neural networks to provide corrective forcing. This approach
includes improvements in the temporal resolution of the simulations.
Learning corrective forcing terms shares similarities with the approach of
DNN-MG, as both methods introduce corrections to the solution through the
right-hand side, which can be seen as corrective forcing terms. However, DNN-MG
goes further by directly correcting the solution in a richer function
space, potentially providing
advantages for goal quantities like drag and lift. The accuracy added by the corrective forcing term on the other hand is limited by the smaller function space.
All these methods have a resemblance to reduced order models (ROMs). An under
resolved simulation, which can be interpreted as a ROM, is augmented by a neural
network. High-fidelity simulations are performed on a small
number of configurations. The high-fidelity data is then utilized to compute a
correction to a coarse/reduced order solution, which yields the enhanced ROM that
accurately approximates solutions for numerous other configurations.
Rather than relying on classical proper orthogonal decomposition
techniques, the idea is to use a simple finite element model that is
refined through the use of a neural network. This results
in a physics-informed, flexible ROM capable of handling various boundary
conditions with enhanced accuracy.
Along these lines,
in~\cite{baigesFiniteElementReducedorder2020} a ROM based on adaptive finite
elements and correction terms combined with artificial neural networks (ANNs) is
presented.
In~\cite{kharaNeuFENetNeuralFinite2021} the authors developed a
mesh based approach to ANN based solvers, leveraging existing FEM theory. The
authors of~\cite{brevisNeuralControlDiscrete2022} propose a neural framework
for optimizing discrete solutions w.\,r.\,t.\ a cost functional in FEM.\@ An ANN
acting as control variable is included in the weak form and the minimization of a cost function yields
desirable approximations where known data or stabilization mechanisms can be
easily incorporated. Ainsworth and Dong~\cite{ainsworthGalerkinNeuralNetworks2021} propose a framework for adaptively
generating basis functions, each represented by a neural network with a single
hidden layer, which may be used in a standard Galerkin method.
In~\cite{mituschHybridFEMNNModels2021} the authors utilize ANNs to represent
unknown physical components of PDEs. For instance, this allows them to recover coefficients
and missing PDE operators from observations. The hybrid model consisting of a
PDE and ANN is then discretized using FEM.\@

The MgNet framework introduced in~\cite{heMgNetUnifiedFramework2019} draws
connections from convolutional neural networks (CNNs) to multigrid methods and thereby improve the design and
understanding of CNNs. The authors of~\cite{luzLearningAlgebraicMultigrid2020}
propose a framework for unsupervised learning of Algebraic Multigrid (AMG)
prolongation operators for linear systems which are defined directly on graphs.
In the multigrid framework, selecting an appropriate smoother is crucial to
ensure convergence and efficiency, but the optimal choice of a smoother is
problem dependent and is in general a major challenge for many applications.
In~\cite{huangLearningOptimalMultigrid2023}, the authors optimize the smoothers
in multigrid methods, which they form by CNNs.

\section{Numerical Methods}
\label{sec:num-meth}
We solve model problems of incompressible viscous flow around a cylinder in 3D domains. This leads to the solution of the Navier-Stokes equations. Our
approach involves a geometric multigrid method with a cell-based Vanka-smoother
as a preconditioner to a Newton-Krylov method. However, any alternative coarse mesh solver can be used as well.

We consider a domain \(\Omega\in \R^d\) with \(d\in \{2,\,3\}\) and
Lipschitz-continuous boundary \(\Gamma\) and a bounded time interval
\([0,\,T]\). For solving the instationary Navier-Stokes equations we seek the
velocity \(\vec v\colon [0,\,T]\times \Omega \to \R^d\) and pressure
\(p\colon [0,\,T]\times \Omega \to \R\), such that
\begin{equation}
  \begin{alignedat}{2}\label{eq:nsstrong}
    \partial_t \vec v + (\vec v\cdot \nabla)\vec v - \nu\Delta \vec v
    +\nabla p &= f \quad &&\text{on } [0,\,T] \times \Omega\\
    \nabla \cdot \vec v &= 0 \quad &&\text{on } [0,\,T] \times \Omega,
  \end{alignedat}
\end{equation}
where \(\nu>0\) is the kinematic viscosity and \(f\) an external force. The
initial and boundary conditions are given by
\begin{equation}\label{eq:boundary}
  \begin{alignedat}{2}
    \vec v(0,\,\cdot) &= \vec v_0(\cdot)\quad &&\text{on }\Omega\\
    \vec v &= \vec v^D \quad &&\text{on } [0,\,T] \times \Gamma^D\\
    \nu(\vec n\cdot\nabla)\vec v - p\vec n &=0 \quad &&\text{in } [0,\,T] \times \Gamma^N,
  \end{alignedat}
\end{equation}
where \(\vec n\) denotes the outward facing unit normal vector on the boundary
\(\partial\Omega\) of the domain. The boundary
\(\Gamma = \Gamma^D \cup \Gamma^N\) is the union of \(\Gamma^D\) with Dirichlet
boundary conditions and \(\Gamma^N\) with Neumann type conditions.

\subsection{Notation and variational formulation of the Navier-Stokes equations}
\label{sec:variational}
By \(L^2(\Omega)\) we denote the space of square integrable functions on the
domain \(\Omega\subset\mathbb{R}^d\) with scalar product \((\cdot \, ,\cdot)\) and
by \(H^1(\Omega)\) those \(L^2(\Omega)\) functions with weak first derivative in
\(L^2(\Omega)\). The function spaces for the velocity and pressure are then
\begin{equation}\label{eq:VQ}
  \begin{alignedat}{1}
     \vec V &\coloneqq \vec v^D + H^1_0(\Omega;\Gamma^D)^d,\quad H_0^1 (\Omega;\Gamma^D)^d\coloneqq \brc{\vec v\in H^1(\Omega)^d\colon \vec v=0 \text{  on } \Gamma^D}\\
    L &\coloneqq \brc{p\in L^2(\Omega),\text{ and, if }\Gamma^N=\emptyset,\; \int_{\Omega}p \drv x = 0},
  \end{alignedat}
\end{equation}
where \(\vec v^D\in H^1(\Omega)^d\) is an extension of the Dirichlet data on
\(\Gamma^D\) into the domain. We normalize the pressure to yield uniqueness, if
Dirichlet data is given on the whole boundary. With these spaces, the
variational formulation of~\eqref{eq:nsstrong} is given by
\begin{equation}
\begin{alignedat}{2}\label{eq:ns}
  (\partial_t \vec v,\,\fvec \phi) + (\vec v\cdot \nabla \vec v,\, \fvec \phi) +
  \nu(\nabla \vec v,\, \nabla \fvec\phi) -
  (p,\,\nabla\cdot \fvec\phi)
  &= (\vec f,\,\fvec \phi)\quad&&\forall \fvec \phi \in H^1_0(\Omega;\Gamma^D)^d\\
  (\nabla \cdot \vec v, \, \xi) &= 0\quad &&\forall \xi \in L\\
  \vec v(0,\,\cdot ) &= \vec v_0(\cdot)\quad &&\text{on }\Omega.
\end{alignedat}
\end{equation}

Let \(\Omega_h\) be a quadrilateral or hexahedral finite element mesh of the
domain \(\Omega\) satisfying the usual requirements on the structural and form
regularity such that the standard interpolation results hold,
compare~\cite[Section 4.2]{richterFluidstructureInteractionsModels2017}.
By \(h_T\) we denote the diameter of an element
\(T\in\Omega_h\) and by \(h\) the maximum diameter of all \(T\in\Omega_h\).

\subsection{Semidiscretization in space}
\label{sec:space-disc}
For the finite element discretization of~\eqref{eq:ns} we choose equal-order continuous finite elements of degree two for velocity and pressure on hexahedral meshes. By \(\mathbb{Q}_h^{(r)}\) we denote the space of continuous functions which are polynomials of maximum directional degree \(r\) on each element \(T\in\Omega_h\).
Then we define the
discrete trial- and test-spaces for the discretization of~\eqref{eq:ns} as
\(\vec v_h,\,\fvec\psi_h \in \vec V_h = [\mathbb{Q}_h^{(2)}]^d\) and
\(p_h,\,\xi_h \in L_h = \mathbb{Q}_h^{(2)}\). Since the resulting equal order
finite element pair \(\vec V_h\times L_h\) does not fulfill the inf-sup
condition, we add stabilization terms of local projection type~\cite{beckerFiniteElementPressure2001}.
We further add convection stabilization, also based on
local projections~\cite{beckerTwolevelStabilizationScheme2004}.
Let $\vec f_h$ and $\vec v^D_h$ be finite element approximations of $\vec f$  and $\vec v^D$ obtained by interpolation in $\vec V_h$.
The
resulting semidiscrete variational problem then reads:
For given data and boundary conditions $\vec f_h,\:\vec v_h^{D}\in
C(I;\,\vec V_{h})$, find $\vec v_h$, $p_h$ such that $\vec v_{h}=\vec v_h^D\:\text{on}\:\bar{I}\times\Gamma_D$ and
\begin{equation}
  \begin{alignedat}{2}\label{eq:disc:ns}
    (\partial_t \vec v_h,\,\fvec\psi_h) + (\vec v_h\cdot \nabla \vec v_h,\, \fvec\psi_h) +
    \nu (\nabla \vec v_h,\, \nabla \fvec\psi_h) -
    (p_h,\,\nabla\cdot \fvec\psi_h)\qquad&\\
    +\sum_{T\in\Omega_h}(\vec v_h\cdot \nabla \vec v_h - \pi_h \vec v_h\cdot
    \nabla \pi_h \vec v_h,\,\delta_T(\vec v_h\cdot \nabla \fvec \psi_h-\pi_h \vec v_h\cdot
    \nabla \pi_h \fvec \psi_h))
    &= (\vec f_h,\,\fvec\psi_h)\quad&&\forall \fvec\psi_h\in  \vec V_h\,,
    \\[5pt]
    (\nabla \cdot \vec v_h, \, \xi_h)
    +\sum_{T\in\Omega_h}\alpha_T (\nabla (p_h-\pi_h p_h),\nabla (\xi_h-\pi_h \xi_h))
    &= 0\quad &&\forall \xi_h \in L_h\,.
  \end{alignedat}
\end{equation}
Here, we denote by $\pi_h:Q^{(2)}_h\to Q^{(1)}_h$ the interpolation into the space of linear finite elements and by $\alpha_T,\delta_T$ two local stabilization parameters specified in~\eqref{stabparams}.

\subsection{Time discretization}
\label{sec:time-disc}
For temporal discretization, the time interval \([0,\,T]\) is split into discrete time steps of uniform size
\[
0=t_0<t_1<\cdots <t_N=T,\: k = t_n-t_{n-1}.
\]
The generalization to a non-equidistant time discretization is straightforward
and only omitted for ease of presentation. We define
\(\vec v_n\coloneqq \vec v_h(t_n)\) and \(p_n\coloneqq p_h(t_n)\) for the fully
discrete approximation of velocity and pressure at time \(t_n\) and apply the
second-order Crank-Nicolson method to~\eqref{eq:disc:ns}, resulting in the
fully discrete problem
\begin{align}
\label{eq:4cranknicholson3}
(\nabla \cdot \vec v_{n}, \, \xi_h)+
\sum_{T\in\Omega_h}\alpha_T (\nabla (p_n-\pi_h p_n),\nabla (\xi_h-\pi_h \xi_h))
&=0 & \forall \xi_h&\in L_h, \nonumber \\[4pt]
\frac{1}{k}(\vec v_n,\,\fvec\phi_h)\,
+{}\frac{1}{2} (\vec v_n\cdot \nabla \vec v_n,\,\fvec\phi_h)
+\frac{\nu}{2}(\nabla \vec v_n,\,\nabla \fvec\phi_h)
-(p_n,\,\nabla \cdot \fvec\phi_h)\quad& \nonumber \\
\qquad \quad \quad+\sum_{T\in\Omega_h}(\vec v_n\cdot \nabla \vec v_n - \pi_h \vec v_n\cdot
    \nabla \pi_h \vec v_n,\,\delta_T(\vec v_n\cdot \nabla \fvec\phi_h-\pi_h \vec v_n\cdot
    \nabla \pi_h \fvec\phi_h)) \nonumber
    &=\frac{1}{k}(\vec v_{n-1},\, \fvec\phi_h)\\
  +\frac{1}{2}(\vec f_n,\fvec\phi_h)
  +\frac{1}{2}(\vec f_{n-1},\,\fvec\phi_h)
  -\frac{1}{2}(\vec v_{n-1}\cdot\nabla \vec v_{n-1},\,\fvec\phi_h)
  &- \frac{\nu}{2}(\nabla \vec v_{n-1},\,\nabla \fvec\phi_h)
  & \forall \fvec\phi_h&\in \vec V_h.
\end{align}
The right hand side only depends on the velocity \(\vec v_{n-1}\) at the last
time step \(n-1\), and we will denote it as \(\vec b_{n-1}\) in the following.

The stabilization parameters $\alpha_T$ and $\delta_T$ depend on the mesh Peclet number and with two parameters $\alpha_0>0$ and $\delta_0\ge 0$ we define them as
\begin{equation}\label{stabparams}
  \alpha_T = \alpha_0 \left(\frac{\nu}{h_T^2}+\frac{\|\vec v_h\|_\infty}{h_T} + \frac{1}{k} \right)^{-1},\quad
  \delta_T = \delta_0 \left(\frac{\nu}{h_T^2}+\frac{\|\vec v_h\|_\infty}{h_T} + \frac{1}{k} \right)^{-1},
\end{equation}
see~\cite{braackStabilizedFiniteElement2007} for details. We set $\alpha_0=0.02$ and $\delta=0.1$.
Introducing the unknown $\vec x_n=(p_n,\,\vec v_n)$, we write Equations~\eqref{eq:4cranknicholson3} more compactly as
\begin{equation}\label{nonlinearshort}
  \begin{aligned}
    {\cal A}_h(\vec x_n) &= \vec F_h,
  \end{aligned}
\end{equation}
where $[{\cal A}_h(\vec x_n)]_i$ and $[\vec F_h]_i$ are the left and right hand sides
of~\eqref{eq:4cranknicholson3} for all test functions
$(\xi_h^i,\,\fvec \phi_h^i)$.

In the presented form, the Crank-Nicolson time discretization is sufficiently robust for smooth
initial data \(v_0\); we refer to~\cite{heywoodFiniteElementApproximation1990} for small
modifications with improved robustness and stability.

\subsection{Solution of the algebraic systems}
\label{sec:algebraic}
The discretization in space and time in Equations~\eqref{eq:4cranknicholson3} leads to a nonlinear system of equations. The nonlinear problem is solved by Newton's method, where we usually omit the dependency on the stabilization parameters in Equations~\eqref{stabparams} on the velocity when setting up the Jacobian. The Jacobian matrix within the Newton iteration is kept for several iterations and also for consecutive time steps and only re-assembled when the convergence rate $\rho= {\lVert{\vec r^{(j+1)}_n}\rVert}/{\lVert{\vec r^{(j)}_n}\rVert}$  of the weak residual $\vec r_n$ in the $j$th Newton iteration deteriorates above a certain threshold $\rho_{\max}$. We usually set $\rho_{\max}$ to $0.1$ or $0.05$. For the linear systems of equations that arise for each Newton iteration, we use a GMRES method~\cite{kelleyIterativeMethodsLinear1995b} that is preconditioned with a geometric multigrid solver~\cite{beckerMultigridTechniquesFinite2000}.

To cope with the structure of the incompressible Navier-Stokes equations and to allow for efficient shared memory parallelization, a Vanka smoother is used within the geometric multigrid. Therefore, small matrices describing the local problems on each mesh element are inverted in a parallel loop with Jacobi coupling. For the quadratic finite elements on hexahedral meshes used in our work, these local matrices are of size $108\times 108$. The complete layout of the algebraic solver and its parallelization is described in~\cite{failerNewtonMultigridFramework2020}.

\section{Neural Networks}
\label{sec:nn}
Various neural network architectures are in principle suitable for our approach. Having previously worked with recurrent neural networks designed for sequence prediction~\cite{margenbergNeuralNetworkMultigrid2022}, we use simpler feedforward neural networks for this work. The simple structure makes them faster to train then the previously used recurrent networks, which require training with backpropagation through time (BPTT) and, consequently, more computationally expensive gradient updates. 
Furthermore, we found that with proper regularization, the feedforward neural network described below performs more consistent than the alternatives we tried. 
Initially, we observed instability when training MLPs with the Adam optimizer. Using the AdamW optimizer~\autocite{loshchilovDecoupledWeightDecay2018}, these issues were resolved and we achieved stability that was on par with recurrent neural networks using backpropagation through time (BPTT). This underscores the importance of regularization in neural network training.

A general feedforward neural network is a composition of multiple parameterized nonlinear functions defined by
\begin{equation}
F(x) \coloneq F(x;\:\mathcal{W}) = (f_L\circ\cdots\circ f_1)(x;\:\mathcal{W})\,,
\end{equation}
where \(\mathcal{W}\) denotes the set of optimizable (learnable) parameters. The \(i\)-th layer of the feedforward network we consider here consists of a weight matrix
\(W_i \in \R^{n_i\times n_{i+1}}\) %
and a nonlinear function \(\sigma\colon \R\to \R\), called the activation function, which is
applied component wise, $\sigma\big(Wx\big)_i = \sigma\big((Wx)_i\big)$.
Furthermore, we employ two techniques common in deep learning to accelerate training and to improve generalization~\cite{zhang2023dive}. Batch normalization $N_i \colon \R^{n_{i+1}} \to \R^{n_{i+1}}$, applied before each activation, re-scales each component by
\[
N_i(x) = \frac{x - \mu_i}{s_i} \odot \gamma_i + \beta_i,
\]
where $x \odot y$ is the Hadamard product of two matrices (or tensors) and
$ \mu_i \in \R^{n_{i+1}}$ and $s_i \in \R^{n_{i+1}} $ are respectively, mean and standard deviation of the samples in a batch estimated during training. The additional parameters $\gamma_i \in \R^{n_{i+1}}$ and $\beta_i \in \R^{n_{i+1}}$ are a learnable, component-wise affine transformation.
In 2D, batch normalization slightly outperformed layer normalization and we therefore also use it in $3D$. See~\cite{Huang2023} for an extensive survey of normalization techniques.

Skip connections are inserted where input and output dimensions match, by adding
the inputs of a layer to its outputs, leading to layers $f_i$ of the structure
\begin{align*}
f_1(x) &= \sigma_1 (N_1( W_1 x)) \\ %
f_i(x) &= \sigma_i (N_i(W_i x)) + x,\quad i = 2,\dots,\,L-1\, \\
f_L(x) &= W_L x.
\end{align*}

\section{The Deep Neural Network Multigrid Method}
\label{sec:dnnmg}
In this section, we review the Deep Neural Network Multigrid (DNN-MG) solver
that we introduced in~\cite{margenbergNeuralNetworkMultigrid2022}. DNN-MG
leverages a deep neural network to predict the correction of a coarse mesh
solution that has been prolongated onto one or multiple finer mesh levels. The
objective is to obtain solutions that are more accurate than solving only on the coarse
mesh, while increasing computational efficiency compared to performing
 a direct simulation on the fine mesh levels. We develop the DNN-MG
solver in a general formulation while maintaining the connection to the
Navier-Stokes equations for which we developed DNN-MG.\@ We begin by providing an
overview of the modifications made in DNN-MG compared to the MG method during a
time step of the Navier-Stokes simulation. We continue by describing the
structure, inputs and outputs of the neural network in DNN-MG.\ In particular,
the network is designed to make local predictions over patches of the mesh.

\LinesNumbered%
\SetFuncSty{textbf} \SetCommentSty{textsf} \SetKwInOut{Input}{Input} %
\SetKwProg{Function}{function}{}{end} %
\SetKwFunction{MG}{Multigrid}%
\SetKwFunction{rhs}{Rhs}%
\begin{algorithm}[ht]
  \Function{DNN-MG}{ %
    \For{all time steps $n$}{%
      \While(\tcp*[f]{Newton-solution of the nonlinear system}){not converged}{
        $\delta z_i =$\label{alg:nkbegin}
        \MG{$L,\,A_{L}^n,\,b_{L}^n,\,\delta z_i$}\tcp*{Solution of the linear
          system (cf.~Section~\ref{sec:algebraic})}%
        $z_{i+1} = z_i + \varepsilon \, \delta z_i$ }%
      \label{alg:nkend}
      \blu{$\tilde{\vec x}_n^{L+J} = \mathcal{P}(\vec x_n^{L}) $}\tcp*{Prolongation on
        level $L+J$}%
      \label{alg:prolong}
      \blu{$\vec r_n^{L+J}=\vec F_h^{L+J}-{\mathcal A}_{h}^{L+J}({\mathcal P}(\vec x_n^{L}))$}\tcp*{Residual
        on level $L+J$ (cf.~\eqref{nonlinearshort})}
      \label{alg:residual}%
      \blu{$\vec d_n^{L+J} =
        \mathcal{L}\circ\mathcal{N}(\mathcal{Q}(\tilde{\vec x}_n^{L+J}),\,\mathcal{Q}
        (\vec r_n^{L+J}),\,\Omega_{h}^{L+J})$}\tcp*{Prediction of velocity
        correction}
      \label{alg:predict}%
      \blu{$b_{n+1}^{L+J} =$
        \rhs{$\tilde{\vec x}_n^{L+J} + \vec d_n^{L+J},f_n,f_{n+1}$}}\tcp*{Set up rhs of~\eqref{eq:4cranknicholson3} for next time step}
      \label{alg:rhs}%
      \blu{$b_{n+1}^{L} = \mathcal{R}(b_{n+1}^{L+J})$}\tcp*{Restriction of rhs
        to level $L$}
      \label{alg:restrictrhs}%
    }%
  }%
  \caption{DNN-MG for the solution of the Navier-Stokes equations. Lines 6-9
    (blue) provide the modifications of the DNN-MG method compared to a classical
    Newton-Krylov simulation with geometric multigrid preconditioning.}\label{alg:dnnmg}
\end{algorithm}
\paragraph{Notation}
In the geometric multigrid, we have a hierarchy of meshes
$\Omega_{h}^{0},\dots,\,\Omega_h^{L+J}$, obtained by successive refinement of the
initial mesh $\Omega_h^0$. We denote the finite element space $\vec V_h$ and
$L_h$ introduced in Section~\ref{sec:space-disc} on level
$l\in\brc{0,\dots,\,L+J}$ by $V_h^l$ and $L_h^l$. With their nested structure,
these spaces reflect the grid hierarchy. Furthermore, let $\Omega_{h}^l$ be the quadrilateral or
hexahedral finite element mesh on level $l\in\brc{0,\dots,\,L+J}$ with
$N_{\text{cells}}^{l}= \operatorname{card}(\Omega_{h}^l)$ cells and the index
set $Z_l$ of all global degrees of freedom
with cardinality $N_{\text{dof}}^{l}\coloneq\operatorname{card}(Z_l)$.
This notation and the description in Section~\ref{sec:dnnnet} is inspired by the description of the Vanka smoother in~\cite{anselmannEfficiencyLocalVanka2023}.

\subsection{Time stepping using DNN-MG}
We present a detailed description of one time step in the simulation of the
Navier-Stokes equations using the DNN-MG solver. The computations involved are
summarized in Algorithm~\ref{alg:dnnmg}.

At the beginning of time step $n$, we solve for the unknown velocity $\vec x_n^{L}$
and pressure $p_n^{L}$ on the coarse level $L$ using the classical Newton-Krylov
simulation (Algorithm~\ref{alg:dnnmg}, lines~\ref{alg:nkbegin}--\ref{alg:nkend})
as described in Section~\ref{sec:algebraic}. Subsequently, we prolongate
$\vec x_n^{L}$ to $\tilde{\vec x}_{n}^{L+J}$ on the finer level $L+J$
(Algorithm~\ref{alg:dnnmg}, line~\ref{alg:prolong}) where a richer function
space $L_h^{L+J}\times \vec V_h^{L+J}$ is available.

Using the prolongated solution, we calculate the residual on level $L+J$ which is part of the
input to the neural network (Algorithm~\ref{alg:dnnmg},
line~\ref{alg:residual}). The residual is calculated according to~\eqref{nonlinearshort}.
We then use the neural network of the DNN-MG solver to
predict the velocity update $\vec d_n^{L+J}$, which represents the difference
$\vec d^{L+J}_{n} = \bar{\vec v}_n^{L+J}- \tilde{\vec v}_n^{L+J}$ between the prolongated
$\tilde{\vec v}_n^{L+J}$ and the unknown ground truth solution $\bar{\vec v}_n^{L+J}$ on
level $L+J$ (Algorithm~\ref{alg:dnnmg}, line~\ref{alg:predict}). The prediction
is based on $\tilde{\vec x}_{n}^{L+J} \in \vec V_h^{L+J}\times L_h^{L+J}$ and $\vec r_n^{L+J}$ and further
utilizes information about the local mesh structure
$\omega\in\R^{n_{\text{Geo}}}$. The $n_{\text{Geo}}$ features of the local geometry are extracted from
the cells in $\Omega_h^{L}$~(cf.~\autocite{margenbergNeuralNetworkMultigrid2022}).

To conclude the $n^{\textrm{th}}$ time step, we compute the right-hand side
$\vec b_{n+1}^{L+J}$ of~\eqref{eq:4cranknicholson3} on level $L+J$ using the
corrected velocity $\tilde{\vec v}_n^{L+J} + \vec d_n^{L+J}$ (Algorithm~\ref{alg:dnnmg},
line~\ref{alg:rhs}) and then restrict it to level $L$
(Algorithm~\ref{alg:dnnmg}, line~\ref{alg:restrictrhs}). This right-hand side,
incorporating the neural network-based correction, becomes part of the solution in the next time step (Algorithm~\ref{alg:dnnmg},
line~\ref{alg:nkbegin}), thereby influencing the overall time evolution of the
flow. This approach of constructing the right-hand side $\vec b_{n+1}^{L+J}$ on level
$L+J$ and subsequently restricting it is a crucial aspect of the DNN-MG solver and 
ensures that the correction on level $L+J$ can have an effect on the time evolution of the flow.

We note that for DNN-MG applied to the Navier-Stokes equations, the pressure is handled
implicitly on level $L$ in the coarse mesh solve and does not receive a
correction by the neural network. However, the pressure is included in the
network's input through the prolongated solution and the residual.
It is indirectly corrected through the
corrections to the velocity.

\subsection{The Neural Network of DNN-MG}\label{sec:dnnnet}
\begin{figure}[tb]
  \includegraphics[width=\linewidth]{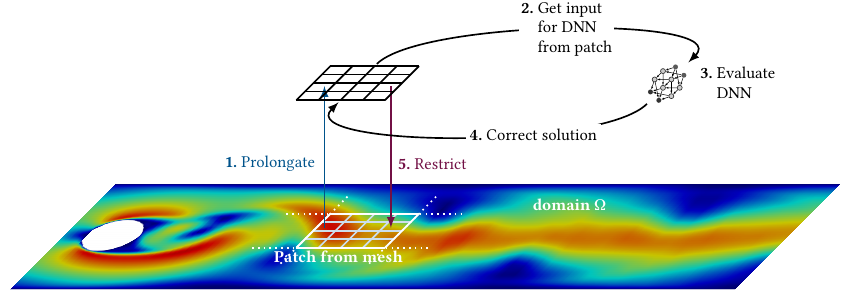}
  \caption{DNN-MG adopts a local approach where the neural network works on
    small neighborhoods of the simulation domain called patches, and
    independently provides a velocity correction for each one.
  } \label{fig:network_patch}
\end{figure}
The neural network component forms the core of the DNN-MG
solver. It plays a crucial role in enhancing computational efficiency,
facilitating fast training procedures, and ensuring robust generalization across
diverse flow regimes beyond the training set. 
We achieve this
through the careful design of a compact neural network architecture with a
moderate number of parameters and a localized, patch-based structure.
Figure~\ref{fig:network_patch} provides an overview of the network component,
illustrating the local approach, which we introduce next.

\paragraph{A patch-based neural network}
To ensure computational efficiency, the DNN-MG approach adopts a neural network
that operates on local patches of the mesh. The network's input and output
are constrained to these patches, which reduces the computational cost compared to a
prediction over the entire domain. The local approach also leads to more compact
neural networks with a smaller number of parameters.

We define a patch $P_{L,\,J}^M\subset\Omega_{h}^{L+J}$ as a collection of cells
$K\in \Omega_{h}^{L+J}$ on level $L+J$. The neural network operates on a single
patch $P_{L,\,J}^M$ and predicts the velocity update
$\vec d_{n,\,P}^{L+J}$ specifically for the degrees of freedom within that
patch. 
Since a single neural network is used in parallel on each patch, all patches must have the same structure and consist of the same number of fine mesh elements.
Formally, let
$L(P_{L,\,J}^M)\times \vec V({P_{L,\,J}^M})\subset L_h^{L+J}\times \vec
V_h^{L+J}$ be the function space associated with the degrees of freedom in patch
$P_{L,\,J}^M$. The neural network learns a mapping
\begin{equation}
  \label{eq:NN}
  \begin{alignedat}{1}
\mathcal{N}: \paran{L(P_{L,\,J}^M)\times \vec V(P_{L,\,J}^M)}^2\times P_{L,J}^M &\to 
    \vec V(P_{L,\,J}^M)\,,\\
    (\tilde{\vec x}_n^{L+J}\restrict{P_{L,\,J}^M},\,\vec r_n^{L+J}\restrict{P_{L,\,J}^M},\,M) &\mapsto \vec d_{n,\,P}^{L+J}\,,
  \end{alignedat}
\end{equation}
that predicts a local velocity update $\vec d_{n,\,P}^{L+J}\in \vec
V(P_{L,\,J}^M)$ based on local data which is also restricted to the patch $P_{L,\,J}^M$.
This prediction is repeated
independently for each patch in the domain, covering the entire domain
$\Omega_h^{L+J}$. The predictions of shared degrees of freedom in adjacent
patches are averaged to ensure consistency. We now formulate the process we just sketched in
more mathematical terms.

We define the patches $P_{L,\,J}^M$ such that each element of mesh level $\Omega_h^L$ will be split into 8 patches. For $J=1$ they consist of one element on level $\Omega_h^{L+1}$  and for $J=2$, $2^d$ elements form one patch.
Additionally, we define
$\Omega_h^{\text{Patch}}\coloneq \brc{P_{L,\,J}^M\;\text{for}\;M\in \Omega_h^{L} }$.
Note that $\operatorname{card}(\Omega_h^{\text{Patch}})=2^d N_{\text{cells}}^L$.
  Combining $2^d$ patches to form the input for the neural network is motivated by a potentially improved approximation quality through a richer set of information. It comes at the cost of larger networks. Finding the optimal balance is subject to future research.
We
denote the set of global degrees of freedom associated with the patch
$P_{L,\,J}^M$ by $Z(P_{L,\,J}^M)\subset Z_{L+J}$. The cardinality of $Z(P_{L,\,J}^M)$ is $N_{\text{dof}}^{P}\coloneq \operatorname{card}(Z(P_{L,\,J}^M))$ and for the $d$-dimensional Navier-Stokes equations with $d$ velocity components and 1 pressure it holds that $N_\text{dof}^P = (d+1)\cdot 3^d$ for $J=1$ and $(d+1)\cdot 5^d$ in the case $J=2$.
Further, the
index set $\hat Z(P_{L,\,J}^M) \coloneq\{0,\dots,\,N_{\text{dof}}^{P}-1\}$ contains all local
degrees of freedom on $P_{L,\,J}^M$. For a given patch $P_{L,\,J}^M$ and a local
degree of freedom with index $\hat\iota \in \hat Z(P_{L,\,J}^M)$ the mapping
\begin{equation}
  \operatorname{dof}:\Omega_h^{L+J}\times \hat Z(P_{L,\,J}^M)  \rightarrow Z_l\,,\quad(P_{L,\,J}^M,\,\hat\iota) \mapsto i \,,
\end{equation}
provides the uniquely defined global index $i \in Z_{L+J}$.

\begin{figure}[t]
  \begin{center}
    \includegraphics[width=\textwidth]{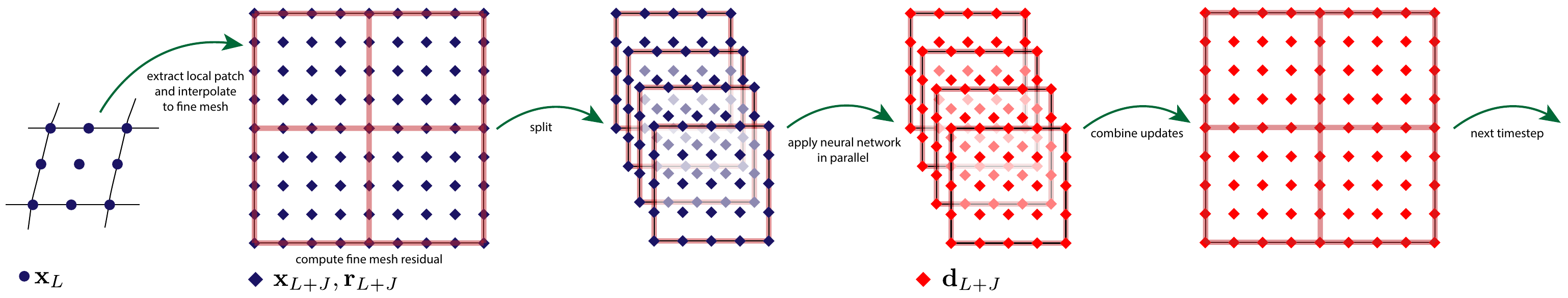}
    \caption{Integration of the neural network for $J=2$: The coarse mesh solution $\mathbf{x}_{L}$ is interpolated to the fine mesh $\mathbf{x}_{L+2}$. Here, the residual $\mathbf{r}_{L+2}$ is computed. The input data $(\mathbf{x}_{L+2},\mathbf{r}_{L+2})$ is split into $2^d$ patches and the network is applied to each patch in parallel to yield the update $\mathbf{d}_{L+2}$. This update is combined on the fine mesh and enters the new time-step.}
    \label{fig:network-sketch}
  \end{center}
\end{figure}

We can now define the prediction of the defect by the neural network in precise
terms. It consists of multiple steps that are illustrated in Fig.~\ref{fig:network-sketch}:
\begin{enumerate}\itemsep1pt \parskip0pt \parsep0pt
\item \emph{Local restriction by}
  $\displaystyle \mathcal Q \colon \R^{N_{\text{dof}}^{L+J}}\to
  \R^{N_{\text{cells}}^L\times N_{\text{dof}}^{P}}$: Using a local restriction
  operator, the globally defined data of size $N_{\text{dof}}^{L+J}$ is
  transformed into local, patchwise defined input of size $N_{\text{dof}}^{P}$
  suitable for the DNN.\@ This is done for each patch, resulting in a matrix of
  size $\R^{N_{\text{cells}}^L\times N_{\text{dof}}^{P}}$.
\item \emph{Patch-wise prediction by the neural network}
  $\displaystyle \mathcal N \colon \R^{(2N_{\text{dof}}^{P} +
    n_{\text{Geo}})}\to \R^{N_{\text{dof}}^{P}}$: For each patch, the neural
  network $\mathcal{N}$ maps the local input data of size
  $(2N_{\text{dof}}^{P} + n_{\text{Geo}})$ to a local output of size
  $N_{\text{dof}}^{P}$, predicting the velocity update for the degrees of
  freedom within the patch $P_{L,,J}^M$. 
  In the output, the degrees of freedom corresponding to the pressure are set to zero.
\item \emph{Global extension by}
  $\displaystyle \mathcal L \colon \R^{N_{\text{cells}}^L \times
    N_{\text{dof}}^{P}} \to \R^{N_{\text{dof}}^{L+J}}$: The locally defined
  output of the DNN is processed by the extension operator $\mathcal{L}$, which
  maps the local data of size $N_{\text{cells}}^L \times N_{\text{dof}}^{P}$ to
  a single globally defined defect vector of size $N_{\text{dof}}^{L+J}$.
\end{enumerate}
In order to define $\mathcal Q$, we need the $P_{L,\,J}^M$-local restriction
$\vec Q\colon\R^{N_{\text{dof}}^{L+J}}\times\Omega_h^{\text{Patch}} \to \R^{N_{\text{dof}}^{P}}$ defined by
\begin{equation}
  \label{Eq:DefR}
  (\vec Q (\vec x,\, P_{L,\,J}^M))[\hat\iota] = \vec x[\operatorname{dof}(P_{L,\,J}^M,\hat\iota)]\,, \quad \text{for} \; \hat\iota \in \hat Z(P_{L,\,J}^M)\,.
\end{equation}
Then, we define $\mathcal Q$ such that it maps a global vector
$\vec x \in \R^{N_{\text{dof}}^{L+J}}$ to a matrix containing local
restrictions to a patch for all $P_0,\dots,\,P_{N_{\text{cells}}^{L}-1}\in\Omega_h^{\text{Patch}}$, stacked on top of each other:
\begin{equation}
  \label{eq:defQ}
  \mathcal Q\colon \R^{N_{\text{dof}}^{L+J}} \to \R^{N_{\text{cells}}^{L}\times N_{\text{dof}}^{P}},\quad \vec x \mapsto \paran{\vec Q_{P}(x,P_0),\dots,\,\vec Q_{P}(x,\,P_{N_{\text{cells}}^{L}-1})}^{\top}
\end{equation}
After applying $\mathcal Q$ to the residual $\vec r_n^{L+J}$ and
$\tilde{\vec x}_n^{L+J}$, we obtain the patch-wise, local data required as input
to the DNN.\@ The network is then evaluated
independently for each patch, with each row of the input representing a
different patch. The output of the DNN then consists of a batch of the patchwise
defect $\vec d_{n,\,P}^{L+J}$. To obtain the global defect vector
$\vec d_{n}^{L+J}$, we use the extension operator $\mathcal{L}$ to transfer this
batch of local data back to the global domain. To this end, we first introduce a
$P_{L,\,J}^M$-local extension operator
$\vec L\colon\R^{N_{\text{dof}}^{P}}\times\Omega_h^{\text{Patch}} \to
\R^{N_{\text{dof}}^{L+J}}$ defined by
\begin{equation*}
  (\vec L( \vec y,\,P_{L,\,J}^M))[\iota]
  =
  \begin{cases}
    \vec y[{\hat \iota}]\,, & \text{if}\; \exists \hat \iota \in \hat Z(P_{L,\,J}^M)\colon \; \iota = \operatorname{dof}(K,\hat \iota),\: K\in P_{L,\,J}^M\,,  \\[1ex]
    0\,, & \text{if}\; \iota \not\in Z_l(P_{L,\,J}^M)
  \end{cases} \,.
\end{equation*}
We further need a scaling vector $\mu\in \R^{N_{\text{dof}}^{L+J}}$ which contains the
reciprocal of the valence of a degree of freedom, i.\,e.\ the number of patches
a degree of freedom is contained in. Then we define $\mathcal L$ such that it
maps $N_{\text{cells}}^{L}$ local vectors to one global vector as
\begin{equation}
  \label{eq:defL}
  \mathcal L \colon \R^{N_{\text{dof}}^{P}\times N_{\text{cells}}^L} \to \R^{N_{\text{dof}}^{L+J}},\quad \vec y\mapsto \mu\odot\sum_{P_{L,\,J}^M\in\Omega_h^{\text{Patch}}} \vec L(\vec y[M],\,P_{L,\,J}^M)\,.
\end{equation}
For the sake of brevity, we used $M$ not only as an element of $\Omega_h^L$
in~\eqref{eq:defL} but also to index the rows of the collection of local
vectors.

This concludes the description of all operations necessary to integrate the
neural network with its local approach into the coarse finite element solver. In particular,
\begin{equation}
  \label{eq:defect-predict}
\mathcal{L}\circ\mathcal{N}(\mathcal{Q}(\tilde{\vec x}_n^{L+J}),\,\mathcal{Q}
        (\vec r_n^{L+J}),\,\Omega_{h}^{L+J})\colon L_h^{L+J} \times \vec
        V_h^{L+J} \to L_h^{L+J} \times \vec V_h^{L+J}\,,\quad
\tilde{\vec x}_n^{L+J}\mapsto \vec d_n^{L+J}\,,
\end{equation}
is fully specified now (cf.~Algorithm~\ref{alg:dnnmg} line~\ref{alg:predict}).
The mapping~\eqref{eq:NN} also fits into this framework and interpreting the
definition~\eqref{Eq:DefR} in a finite element space context, the restrictions
$\tilde{\vec x}_n^{L+J}\restrict{P_{L,\,J}^M}$ and
$\vec r_n^{L+J}\restrict{P_{L,\,J}^M}$ are obtained. 

The neural network operates on
each patch independently, which allow us to take advantage of parallelization
through the localized nature of the problem. According to~\eqref{eq:defQ}, the local data is stacked along the
first dimension to form a batch as processed by the neural network.
Therefore, evaluating both the extension operator $\mathcal Q$ and the network $\mathcal{N}$ can be done in embarrassingly parallel fashion.
The local approach also allows DNN-MG to handle a wide range
of flow scenarios while maintaining computational efficiency, as we will show
later in Section~\ref{sec:num-example}. The rest of this section consists of a
detailed description of the inputs to the neural network and the training
methodology.

\paragraph{Neural network inputs}
The effectiveness of the neural network prediction hinges upon the selection
of appropriate network inputs that provide enough information
regarding the coarse flow characteristics at level $L$ and the local patch
geometry, akin to the significance of mesh cell geometry in classical finite
element simulations. A well-informed choice of these inputs holds the potential to
keep the number of parameters in the neural network small, so that also the
evaluation time during runtime is small and the data and time prerequisites
for training are minimized.
Based on our previous work~\cite{margenbergNeuralNetworkMultigrid2022}, we use the following inputs to the neural network:
\begin{itemize}\itemsep1pt \parskip0pt \parsep0pt
\item nonlinear residuals
  $\vec r_n^{L+J}\restrict{P_{L,\,J}^M}=\paran{\vec F_h^{L+J}-{\mathcal A}_h^{L+J}(\vec{\tilde x}_n^{L+J})}\restrict{P_{L,\,J}^M}\in \R^{N_{\text{dof}}^{P}}$ of
  the prolongated coarse mesh solution for~\eqref{eq:4cranknicholson3};
\item the velocities and pressure
  $\vec{\tilde x}_n^{L+J}\restrict{P_{L,\,J}^M} \in \R^{N_{\text{dof}}^{P}}$ on mesh level $L+J$;
\item the geometry of the cells, which for the 3D case consist of
  \begin{itemize}\itemsep1pt \parskip0pt \parsep0pt
  \item the edge lengths $h^c\in \R^{12}$;
  \item the lengths of the diagonals connecting the vertices which do not share
    the same face $a^{c} \in \R^{4}$;
  \item for each vertex, the average of the 3 angles between the faces
    $\alpha^c\in\R^8$.
  \end{itemize}
\end{itemize}
The nonlinear residual plays a crucial role as it quantifies the local error or
defect at each degree of freedom on level $L+J$ and therefore carries essential
information about the underlying partial differential equation. Furthermore, the
velocity and pressure fields contain details about the current solution and
provide knowledge about the flow within the current patch. It is important to
note that the geometric information of the cells also plays a significant role in
our approach, as we impose no specific restrictions on the cell geometries, except for the
standard shape regularity requirements commonly employed in finite element
methods. The neural network hence has to make predictions that are adapted to the specific cell.

\paragraph{Training of DNN-MG}
The training of DNN-MG is based on simulations of the Navier-Stokes equations
for which a multigrid representation of the velocity $\vec v_n^l$ with two levels $L$
and $L+J$ is available. The velocity $\vec v_n^{L+J}$ thereby serves as ground truth
$\bar{\vec v}_n^{L+J}$. The goal of the training is to optimize the network
parameters such that the $\ell^2$-norm
$\big\Vert (\tilde{\vec v}_n^{L+J} + \vec d_n^{L+J}) - \bar{\vec v}_n^{L +J } \big\Vert_2$
of the difference between the predicted velocity
$\tilde{\vec v}_n^{L+J} + \vec d_n^{L+J}$ (i.\,e.\ the velocity after the correction)
and the ground truth $\bar{\vec v}_n^{L}$ is minimized throughout the simulation.
Typically, the training data for the network comprises only a few snapshots of
simulation data from simple scenarios, as we will see in the next section. 
This is sufficient thanks to the local approach where even a few global snapshots 
provide a rich set of different local flow scenarios.

\section{Numerical examples}
\label{sec:num-example}

In this section, we will document different numerical
simulations. The test cases are based on a well established 3D
Navier-Stokes benchmark
problem~\cite{schferBenchmarkComputationsLaminar1996,braackSolutions3DNavierStokes2006}. While
the neural network is trained on the original benchmark configuration,
we study its performance and generalization on a set of modified
problems. As perturbations, we consider an obstacle with an elliptical cross section and increasing or
decreasing the Reynolds number by varying the viscosity (cf.~Figure~\ref{fig:3d-1-geo}).
We further consider a substantial change to the problem geometry by introducing a second obstacle.
There, we again increase and decrease the Reynolds number (cf.~Figure~\ref{fig:2obs-geo}).
It is important to stress that the network is trained only once on the original benchmark configuration and no retraining is applied for the further test cases.

\subsection{Implementation aspects}
\label{sec:implementation}
The DNN-MG method is implemented using two libraries: The FEM library
Gascoigne3D~\cite{beckerFiniteElementToolkit} and the machine learning library PyTorch~\cite{paszkePyTorchImperativeStyle2019}.
The Newton-Krylov geometric multigrid
implemented in Gascoigne3D method has been shown to be efficient in~\cite{failerNewtonMultigridFramework2020}.

We implement the deep neural networks with \texttt{libtorch}, the \texttt{C++} interface of
\texttt{PyTorch}. Using \texttt{C++} for the
implementation allows us to couple the numerical methods and neural networks in
a performant manner without unnecessary overhead. PyTorch also supports
parallelization with MPI, which is used here to distribute training on multiple
GPUs.

\subsection{The 3D Benchmark Flow Around a Cylinder}
\label{sec:3d-benchmark}

\begin{figure}[t]
  \includegraphics{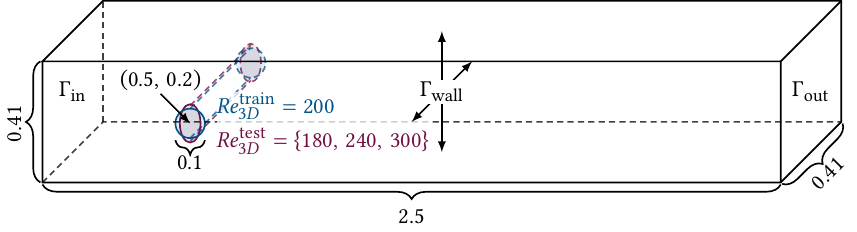}
  \caption{Geometry of the training and test scenario with a parabolic inflow
    profile at $\Gamma_{\textrm{in}}$, do-nothing boundary conditions at the
    outflow boundary $\Gamma_{\textrm{out}}$ and no-slip conditions on the walls
    $\Gamma_{\textrm{wall}}$. The center of the obstacle is at $(0.5,\,0.2)$.
  }\label{fig:3d-1-geo}
\end{figure}

We examine a variant of the three-dimensional flow benchmark presented in~\cite{schferBenchmarkComputationsLaminar1996} in a setting similar to the
\emph{3D-2Z} one, with different Reynolds numbers \(Re \in \{180,\,200,\,240,\,300\}\). We
show the geometry in Figure~\ref{fig:3d-1-geo}.
The meshes used in this paper are generated using \texttt{deal.ii} functions for
creating triangulations of basic
geometries~\cite{arndtDealIILibrary2023}. We generate unstructured, pre-adapted meshes.
While the original benchmark description considers an obstacle with circular cross-section, we allow
elliptical ones for testing the generalizability of the approach. We define the
Reynolds number as
\begin{equation}\label{reynolds}
  \text{Re} = \frac{\bar{\vec v}\cdot L}{\nu},
\end{equation}
where $\bar{\vec v}$ is the mean inflow velocity and $L$ the length of the major
axis of the obstacle. When the cross-section of the cylinder is circular, $L$
represents its diameter. If the cross-section is elliptical,
$L$ corresponds to the height of the obstacle, as the major axis is always
parallel to the $y$-axis.

The time step size is chosen as \(k=\num{0.008}\) on the
interval \(I=(0,8]\). The initial velocity is \(\vec{v}(0) = \vec{0}\), and the
inflow boundary condition is prescribed for \(t>0\), with a smooth startup process.
The flow is driven by a Dirichlet
profile \(\vec v=\vec v^D\) at the left boundary \(\Gamma_{\text{in}}\) given by
\begin{equation}
  \label{eq:inflow}
  \begin{aligned}
    \vec v_{in}^{3d}(t,\,y,\,z)&=\omega(t)\frac{y(H-y)(H^2-z^2)}{(H/2)^2H^2}\frac{9}{8}\bar{\vec v}_{3d} \ \text{  on  } \ \Gamma_{in}\coloneqq
                              0\times [0,\,H],
    \\[4pt]
    \; \omega(t) &= \begin{cases} \frac{1}{2}- \frac{1}{2}\cos( 5\pi t), & t\le \frac{1}{5} \\ 1 & t>\frac{1}{5}.
                    \end{cases}
  \end{aligned}
\end{equation}
\(H=0.41\) is the height and width of the channel and \(\omega(t)\) acts
as a regularization during the startup phase. On the wall boundary \(\Gamma_{wall}\)
and on the obstacle we use no-slip boundary conditions. On the outflow boundary
\(\Gamma_{\text{out}}\) we use a do-nothing outflow condition~\cite{heywoodArtificialBoundariesFlux1992}.
The mean flow rate is \(\bar{\vec v}_{3d}=\frac{1}{\abs{\Gamma_{\text{out}}}}\int_{\Gamma_{\text{out}}}\vec v\drv s= 1\) and the viscosity is \(\nu = \num{5.e-4}\). The training data is generated considering the circular obstacle, i.\,e.\ $L=0.1$ which results in the Reynolds number $\text{Re}=200$. For testing, the elliptical cross-section with $L=0.12$ is considered, leading to the slightly increased Reynolds number $\text{Re}=240$. See Section~\ref{sec:flow-pred} for details.

As in the classical benchmark, we analyze the drag and lift coefficients of the obstacle
\begin{equation}\label{functionals}
  J_\text{drag} = \int_\Gamma
  \big( \nu\nabla\vec v\vec n-p\vec n\big)\vec e_x\,\text{d}s,\quad
  J_\text{lift} = \int_\Gamma
  \big( \nu\nabla\vec v\vec n-p\vec n\big)\vec e_y\,\text{d}s.
\end{equation}
To compute the drag and lift coefficients, we use the Babuszka-Miller trick;
cf.~\cite{braackSolutions3DNavierStokes2006,babukaPostprocessingApproachFinite1984}
and rewrite the surface integrals over the volume for obtaining super-convergence of fourth order.

\begin{table}[t]
  \caption{Spatial meshes with increasing number of
    multigrid levels and degrees of freedom (DoF). For the DNN-MG setups we
    indicate the size of the patches (counted as number of fine mesh elements)
    where the network is applied.}
    \label{tab:org62ab3c1}
  \centering
  \begin{tabular}{lcrr}
    \toprule
    & Level & \# DoF& Patch size\\[0pt]
    \midrule
    Coarse & MG(\(3\)) & \(\num{390720}\) & \num{0}\\[0pt]
    Reference for $\text{DNN-MG}(3+1)$ & MG(\(4\)) & \(\num{3003520}\) &\num{8}\\[0pt]
    Reference for $\text{DNN-MG}(3+2)$ & MG(\(5\)) & \(\num{23546112}\)& \num{64}\\[0pt]
    \bottomrule
  \end{tabular}
\end{table}

\subsection{Neural Network parameters}
\label{sec:nn-parameter}
DNN-MG uses a neural network that operates on patches. In our experiments, a
patch is once or twice refined cell on level \(L\). In
Table~\ref{tab:orgef2be82} we list the dimensions and number of trainable
parameters of the resulting neural networks.
In determining a suitable architecture for our neural networks, we found that
hidden layers with a size of 512 or 750 strike a balance between performance
and model complexity. In our experiments, the networks with a width of 750
performed best, and we include networks with a hidden layer size of 512 only for
comparsion. We generally use 8 hidden layers of the same size.
Smaller networks were more difficult to train and performed less
  consistent than the ones described above.
Our experiments consistently supported the benefits of deeper networks, which we
consider a crucial factor in accurately predicting various flow regimes.
While more extensive tuning might lead to smaller networks, the small overhead
of network evaluation does not justify the required effort.
The hyperparameters were determined based on our prior exeperiments
in 2D (cf.~\cite{margenbergDeepNeuralNetworks2021}). As they proved to be effective, we
omit an extensive hyperparameter analysis. However, this is planned for a
manuscript currently in preparation.

\begin{table}[t]
  \caption{The input size \(N_{\text{in}}\), hidden size
    \(N_{\text{hidden}}\), output size \(N_{\text{out}}\) and the number of
    trainable parameters of the feedforward network, for different numbers of
    predicted levels.}
    \label{tab:orgef2be82}
  \centering
  \begin{tabular}{rrrrr}
    \toprule
    Predicted levels & \(N_{\text{in}}\) & \(N_{\text{hidden}}\) & \(N_{\text{out}}\) & Trainable parameters\\[0pt]
    \midrule
    1 & 240 & 512 & 81 & \num{2007552}\\[0pt]
    2 & 1024 & 512 & 375 & \num{2559488}\\[0pt]
    1 & 240 & 750 & 81 & \num{4190250}\\[0pt]
    2 & 1024 & 750 & 375 & \num{4998750}\\[0pt]
    \bottomrule
  \end{tabular}
\end{table}

\subsection{Neural Network training}
\label{sec:nn-train}

Training data are generated on grids with three successive levels of refinement, see Table~\ref{tab:org62ab3c1}. Here, we denote by MG(\(3\)) the coarse grid. The grid levels MG(\(4\)) and MG(\(5\)) are the training data for a network prediction over one, respectively, over two grid levels.

A single simulation is sufficient as training data set since, 
by the patch-wise application of the network, it provides
\(N_c \times N_T\) training items, where \(N_c\) is the number of patches and
\(N_T\) is the length of the time series. From the experience gained in the 2D
case~\cite{margenbergNeuralNetworkMultigrid2022}, we consider a small subinterval of
\(I=[0,\,8]\). Here we use \(t \in [4,\,7]\) resulting in \(N_T=375\). We
consider the whole spatial domain for training, which results in $N_c=12210$
patches. The interval from \([2,\,4]\) is used as validation set. Therefore, the
training dataset contains $4.57875\times 10^6$ samples and the validation set
consists of $3.0525\times 10^6$ samples. The full dataset encompasses 1\,TB of
data. Table~\ref{tab:stats} show some
descriptive statistics of the velocity and pressure in the training and
validation dataset. In the Appendix~\ref{app:stat} we include the same
statistics for the test cases. The data on which the network is tested on
is significantly different from the training data (cf.~Table~\ref{tab:stats1}).

\begin{table}[t]
  \caption{Statistics of the training dataset on the left and the validation
    dataset on the right.}\label{tab:stats}
  \small
  \sisetup{round-mode=places,round-precision=2,exponent-mode=scientific,print-zero-exponent=false,tight-spacing=true}
  \centering
  \begin{minipage}{0.49\textwidth}
    \begin{center}
      \begin{tabular}{crrrr}
        \toprule
        & max & mean & min & $\sigma$ \\[0pt]
        \midrule
        $\vec v_{x}$ & \num{2.609} & \num{0.978}     & -\num{0.7247}& \num{0.7074} \\[0pt]
        $\vec v_{y}$ & \num{1.663} & \num{0.0005974}  & -\num{1.676} & \num{0.2583} \\[0pt]
        $\vec v_{z}$ & \num{1.245} & \num{2.433e-11}  & -\num{1.245} & \num{0.1705} \\[0pt]
        $p$          & \num{3.397} & \num{0.2389}     & -\num{1.303} & \num{0.5443} \\[0pt]
        \bottomrule
      \end{tabular}
    \end{center}
  \end{minipage}
  \begin{minipage}{0.49\textwidth}
    \begin{center}
      \begin{tabular}{crrrr}
        \toprule
        & max & mean & min & $\sigma$ \\[0pt]
        \midrule
        $\vec v_{x}$ & \num{2.558} & \num{0.978}     & -\num{0.5064}  & \num{0.7003} \\[0pt]
        $\vec v_{y}$ & \num{1.606} & \num{0.0008537} & -\num{1.621}   & \num{0.2479} \\[0pt]
        $\vec v_{z}$ & \num{1.024} & \num{4.912e-13} & -\num{1.024}  &  \num{0.1513} \\[0pt]
        $p$          & \num{3.605} & \num{0.3772}    & -\num{1.182}  & \num{0.5733} \\[0pt]
        \bottomrule
      \end{tabular}
    \end{center}
  \end{minipage}
\end{table}

The networks were trained using 2 GPU nodes, each with 2 Nvidia A100 GPUs and 2
Intel Xeon Platinum 8360Y CPUs. The number of MPI processes were chosen equal to
the number of GPUs such that one MPI process uses one GPU and CPU.\@ Using this
setup, the training of the feedforward network takes one day, regardless of the
network size. Although the data loading process is implemented efficiently and
the average load of the GPUs ranges from $80\%$ for small networks to $90\%$ for
large networks, the limiting factor in terms of performance is the memory
bandwidth.

\begin{figure}[t]
  \centering \includegraphics{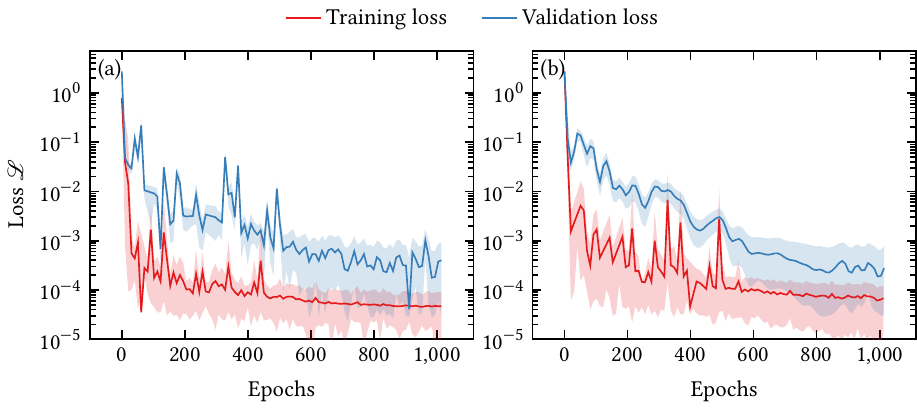}
  \caption{The average training and validation loss of the feedforward network
    architecture with hidden size (a) 512 and (b) 750, plotted over the number of epochs. Due to the data
    efficiency of DNN-MG and the large number of batches
    ($\approx \num{900000}$) relative to the network size, we were able to
    achieve a good training loss after a single epoch and there is only little
    improvement in the training loss beyond that point.}\label{fig:so-losses}
\end{figure}

As loss function \(\mathcal{L}\) we employed a simple \(\ell^2\)-loss
\[
  \mathcal{L} = \sum_{n=1}^{N_T} \frac{1}{N_c} \sum_{c=1}^{N_c}
  l^2(\tilde{v}_{n,c}^{L+J} + d_{n,c}^{L+J} ,\,\bar {v}_{n,c}^{L+J})^2 =
  \sum_{i=1}^{N_T} \frac{1}{N_c} \sum_{c=1}^{N_c} \norm{ \big(
    \tilde{v}_{n,c}^{L+J} + d_{n,c}^{L+J} \big) - \bar{v}_{n,c}^{L+J}
  }_2^2,\quad J\in \{1,\,2\} .
\]
The parameter $J$ is the number of levels we skip with the prediction by the DNN.\@
We use a Tikhonov-regularization with a scaling factor \(\alpha=10^{-5}\). To
optimize our model, we utilize the AdamW~\cite{loshchilovDecoupledWeightDecay2018} optimizer with a
maximum of 1000 epochs. Finally, we select the model with the lowest validation
loss. The convergence of the training for the different network
configurations is reported in Figure~\ref{fig:so-losses}.

\subsection{Flow prediction}
\label{sec:flow-pred}
For testing, we use the flow around an elliptic obstacle with an increased height
of \(0.12\) compared to the training data (see Figure~\ref{fig:3d-1-geo}). The resulting finite element meshes have the same structure as in the training case but the elements are distorted. We
observe that DNN-MG is indeed able to predict high frequency fluctuations that
are not visible in the coarse mesh solution. In particular in the vicinity of
the obstacle the quality of the solution is strongly enhanced with distinct
features in the wake being apparent in the DNN-MG simulation.

\begin{figure}
  \includegraphics{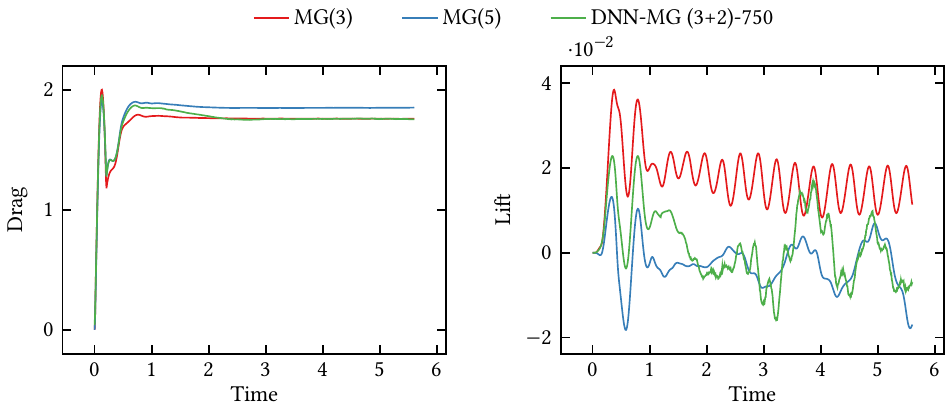}
  \caption{The drag (left) and lift (right) functionals for the channel with one
    obstacle with an elliptical cross section at $\text{Re}=180$ are shown. The
    results are presented for three discretizations: The coarse mesh
    $\text{MG}(3)$, the reference $\text{MG}(5)$, and the $\text{DNN-MG}(3+2)$
    improved by a neural network.}\label{fig:dl-3d-1-re180}
  \captionof{table}{The $\min$, $\max$, mean and amplitude of drag and lift
    functionals measured in the interval $[2,8]$ in the same setting as above in Figure~\ref{fig:dl-3d-1-re180}. Amplitude is defined as $\text{amp}=\max-\min$.
    The best results in comparison to the $MG(5)$ solution are highlighted with
    color coding: \bestf{Reference}, \bests{1st best}, \bestt{2nd
      best}.}\label{tab:dl-3d-1-re180}
  \centering
  \begin{tabular}{lrrrrrrrr}
    \toprule
    Type & \(\min C_{d}\) & \(\max C_{d}\) & \(\widebar{C_{d}}\) & amp \(C_{d}\) & \(\min C_{l}\) & \(\max C_{l}\) & \(\widebar{C_{l}}\) & amp \(C_{l}\)\\[0pt]
    \midrule
    $\text{MG}(3)$ & 1.7571 & 1.7662 & 1.7591 & \bestt{0.0091} & 0.0083 & 0.0235 & 0.0157 & 0.0152\\[0pt]
    $\text{MG}(4)$ & 1.7606 & 1.7713 & 1.7638 & \bests{0.0107} & -0.0016 & 0.0208 & 0.0104 & \bestt{0.0224}\\[0pt]
    $\text{MG}(5)$ & \bestf{1.8472} & \bestf{1.8595} & \bestf{1.8494} & \bestf{0.0123} & \bestf{-0.0178} & \bestf{0.0069} & \bestf{-0.0033} & \bestf{0.0247}\\[0pt]
    $\text{DNN-MG}(3+1)$-512 & 1.7082 & 1.7575 & 1.7186 & 0.0493 & 0.0025 & 0.0241 & 0.0136 & 0.0216\\[0pt]
    $\text{DNN-MG}(3+2)$-512 & \bests{1.8593} & \bests{1.9112} & \bests{1.8742} & 0.0519 & \bests{-0.0223} & 0.0025 & \bestt{-0.0119} & \bests{0.0248}\\[0pt]
    $\text{DNN-MG}(3+1)$-750 & 1.7257 & 1.7652 & 1.7374 & 0.0394 & -0.0235 & \bestt{0.0045} & -0.0049 & 0.0280\\[0pt]
    $\text{DNN-MG}(3+2)$-750 & \bestt{1.7657} & \bestt{1.8052} & \bestt{1.7774} & 0.0394 & \bestt{-0.0235} & \bests{0.0045} & \bests{-0.0049} & 0.0280\\[0pt]
    \bottomrule
  \end{tabular}
\end{figure}

\begin{figure}
  \centering%
  \includegraphics{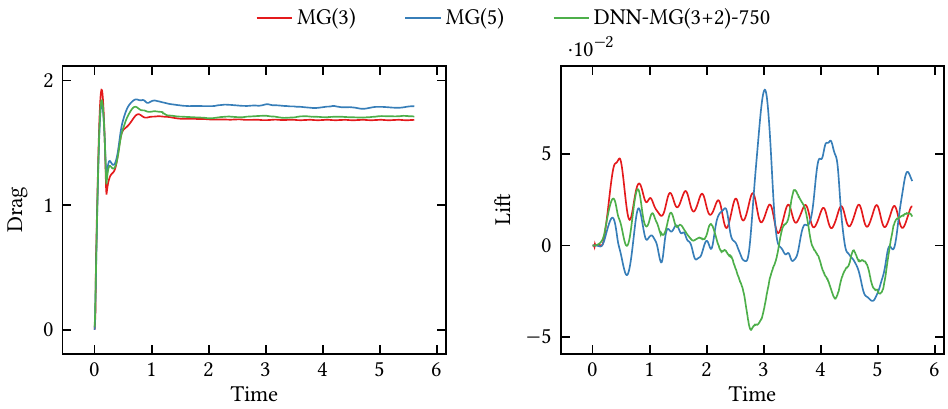}
  \caption{The drag (left) and lift (right) functionals for the channel with one
    obstacle with an elliptical cross section at $\text{Re}=240$ are shown. The
    results are presented for three discretizations: The coarse mesh
    $\text{MG}(3)$, the reference $\text{MG}(5)$, and the $\text{DNN-MG}(3+2)$
    improved by a neural network.}\label{fig:dl-3d-1-re240}
  \captionof{table}{The $\min$, $\max$, mean and amplitude of drag and lift
    functionals in the same setting as above in Figure~\ref{fig:dl-3d-1-re240}.
    The best results in comparison to the $MG(5)$ solution are highlighted with
    color coding: \bestf{Reference}, \bests{1st best}, \bestt{2nd
      best}.}\label{tab:dl-3d-1-re240}
  \begin{tabular}{lrrrrrrrr}
    \toprule
    Type & \(\min C_{d}\) & \(\max C_{d}\) & \(\widebar{C_{d}}\) & amp \(C_{d}\) & \(\min C_{l}\) & \(\max C_{l}\) & \(\widebar{C_{l}}\) & amp \(C_{l}\)\\[0pt]
    \midrule
    $\text{MG}(3)$ & 1.6810 & 1.6920 & 1.6840 & 0.0111 & 0.0081 & \bests{0.0495} & 0.0255 & 0.0415\\[0pt]
    $\text{MG}(4)$ & 1.6783 & 1.6938 & 1.6869 & 0.0155 & \bestt{-0.0204} & 0.0373 & \bestt{0.0098} & 0.0576\\[0pt]
    $\text{MG}(5)$ & \bestf{1.7724} & \bestf{1.8096} & \bestf{1.7910} & \bestf{0.0372} & \bestf{-0.0303} & \bestf{0.0852} & \bestf{0.0123} & \bestf{0.1155}\\[0pt]
    $\text{DNN-MG}(3+1)$-512 & 1.6305 & 1.6748 & 1.6425 & \bests{0.0444} & -0.0364 & 0.0224 & -0.0002 & 0.0588\\[0pt]
    $\text{DNN-MG}(3+2)$-512 & \bestt{1.6866} & \bests{1.7931} & \bests{1.7537} & 0.1066 & -0.0347 & 0.0363 & 0.0089 & \bestt{0.0707}\\[0pt]
    $\text{DNN-MG}(3+1)$-750 & 1.6215 & 1.6664 & 1.6309 & 0.0448 & 0.0008 & 0.0304 & \bests{0.0146} & 0.0296\\[0pt]
    $\text{DNN-MG}(3+2)$-750 & \bests{1.6966} & \bestt{1.7167} & \bestt{1.7079} & \bestt{0.0200} & \bests{-0.0345} & \bestt{0.0426} & 0.0059 & \bests{0.0774}\\[0pt]
    \bottomrule
  \end{tabular}
\end{figure}

\begin{figure}
  \centering%
  \includegraphics{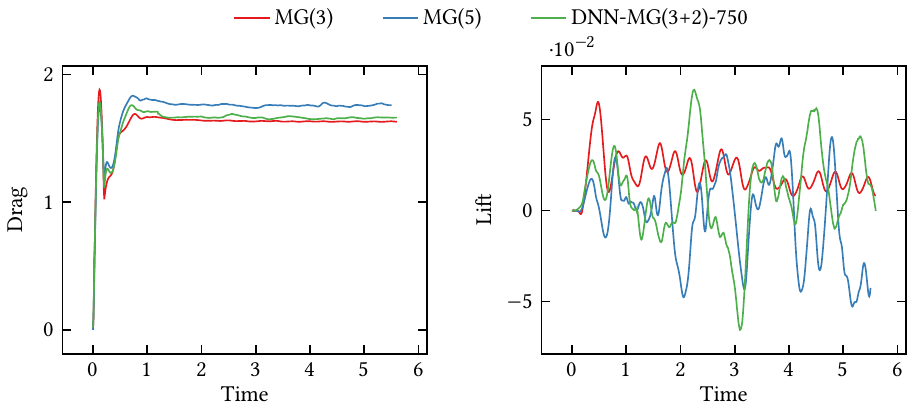}
  \caption{The drag (left) and lift (right) functionals for the channel with one
    obstacle with an elliptical cross section at $\text{Re}=300$ are shown. The
    results are presented for three discretizations: The coarse mesh
    $\text{MG}(3)$, the reference $\text{MG}(5)$, and the $\text{DNN-MG}(3+2)$
    improved by a neural network.}\label{fig:dl-3d-1-re300}
  \captionof{table}{The $\min$, $\max$, mean and amplitude of drag and lift
    functionals in the same setting as above in Figure~\ref{fig:dl-3d-1-re300}.
    The best results in comparison to the $MG(5)$ solution are highlighted with
    color coding: \bestf{Reference}, \bests{1st best}, \bestt{2nd
      best}.}\label{tab:dl-3d-1-re300}
  \begin{tabular}{lrrrrrrrr}
    \toprule
    Type & \(\min C_{d}\) & \(\max C_{d}\) & \(\widebar{C_{d}}\) & amp \(C_{d}\) & \(\min C_{l}\) & \(\max C_{l}\) & \(\widebar{C_{l}}\) & amp \(C_{l}\)\\[0pt]
    \midrule
    $\text{MG}(3)$ & \bestt{1.6291} & 1.6449 & \bestt{1.6339} & 0.0159 & 0.0081 & 0.0337 & 0.0185 & 0.0256\\[0pt]
    $\text{MG}(4)$ & 1.6125 & 1.6507 & 1.6336 & \bestt{0.0381} & \bests{-0.0592} & 0.0724 & \bestt{0.0142} & 0.1316\\[0pt]
    $\text{MG}(5)$ & \bestf{1.7375} & \bestf{1.7794} & \bestf{1.7578} & \bestf{0.0419} & \bestf{-0.0527} & \bestf{0.0404} & \bestf{-0.0058} & \bestf{0.0931}\\[0pt]
    $\text{DNN-MG}(3+1)$-512 & 1.4532 & 1.5883 & 1.4985 & 0.1351 & -0.0047 & \bests{0.0355} & 0.0161 & 0.0401\\[0pt]
    $\text{DNN-MG}(3+2)$-512 & 1.4923 & \bests{1.7125} & 1.5776 & 0.2202 & -0.0071 & \bestt{0.0656} & 0.0266 & \bestt{0.0728}\\[0pt]
    $\text{DNN-MG}(3+1)$-750 & 1.5660 & 1.6347 & 1.5858 & 0.0687 & -0.0235 & 0.0714 & 0.0184 & \bests{0.0950}\\[0pt]
    $\text{DNN-MG}(3+2)$-750 & \bests{1.6490} & \bests{1.6910} & \bests{1.6626} & \bests{0.0419} & \bestt{-0.0655} & 0.0664 & \bests{0.0109} & 0.1319\\[0pt]
    \bottomrule
  \end{tabular}
\end{figure}

\begin{figure}
  \centering \includegraphics{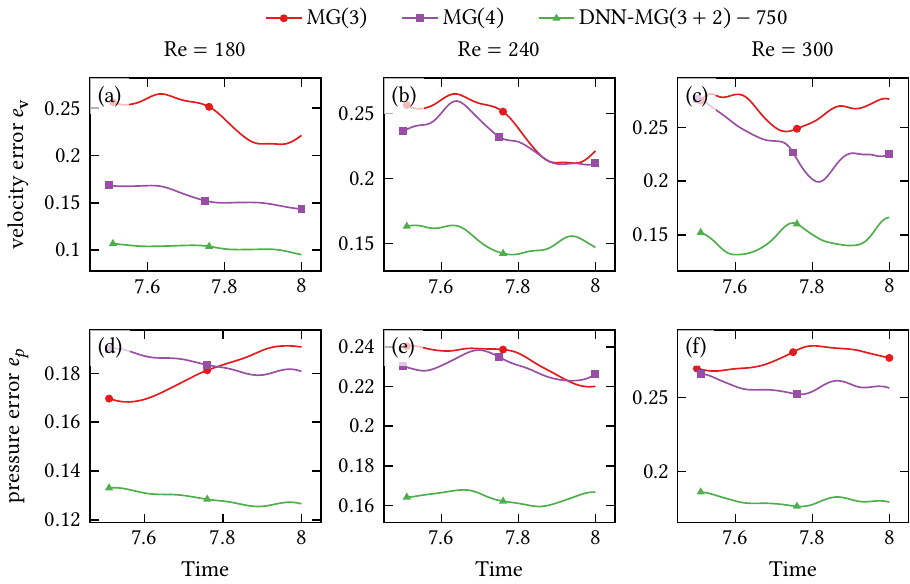}
  \caption{Relative velocity errors (first row) and pressure errors (second row)
    for coarse solutions with and without ANN correction ($\text{MG}(3)$ and
    $\text{DNN-MG}(3+2)-750$) compared to the reference solution on level $L+2$
    for different Reynolds numbers. As further evidence of the innovative
    capacity of DNN-MG we show the errors of the $\text{MG}(4)$ solution as
    well, which we consistently outperform.}\label{fig:vp-errors-1obs}
  \vspace{.5em}
  \includegraphics{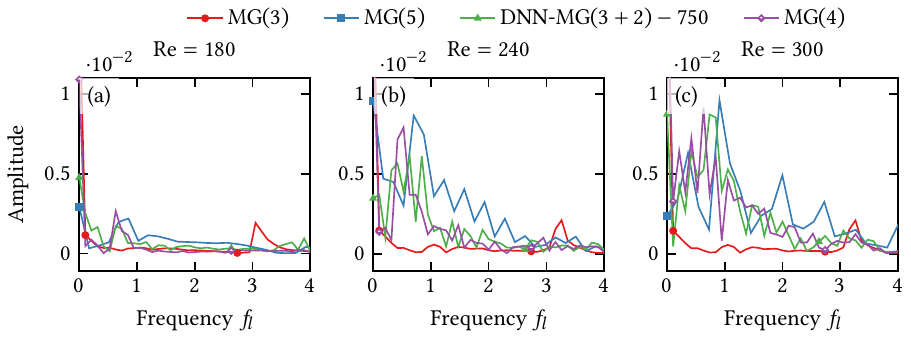}
  \caption{The spectrum of the lift functional for the channel flow with one
    obstacle at different Reynolds numbers. We plot the results of the coarse
    mesh $\text{MG}(3)$, the reference $\text{MG}(5)$, the intermediate level
    $\text{MG}(4)$ and the $\text{DNN-MG}(3+2)-750$ improved by a neural
    network.}\label{fig:lspec1obs}
\end{figure}

In addition to a viscosity of \(\nu=\num{5.e-4}\) and \(Re=240\), we
consider fluids with \(\nu=\num{4e-4}\) corresponding to a Reynolds number of
\(Re=300\) and \(\nu=\num{6.67e-4}\) which corresponds to a Reynolds number of
180. Both of these Reynolds numbers are calculated for the elliptic obstacle
with a height of \(0.12\), which is our typical test case. The network
is trained only for \(Re=200\) using the obstacle with circular cross-section.
\begin{description}[style=unboxed,leftmargin=0cm]\itemsep1pt \parskip0pt \parsep0pt
  \item[Reynolds Number $\boldsymbol{180}$.] 
As shown in Figure~\ref{fig:dl-3d-1-re180} and Table~\ref{tab:dl-3d-1-re180}, 
  DNN-MG is able to predict the dynamic behavior of flow at $\text{Re}=180$,
slightly enhancing the drag coefficient and significantly improving the lift
compared to coarse solutions.
For finite element simulations, the lift is the functional that is usually significantly more difficult to approximate~\cite{braackSolutions3DNavierStokes2006}. The overall flow dynamics are
reconstructed effectively, albeit with temporal deviations. We note that it is
imperative to use DNN-MG with two level predictions, as single level predictions
yield an unreliable velocity and pressure field although drag and lift
functionals in Table~\ref{tab:dl-3d-1-re180} may suggest otherwise.
We refer to Section~\ref{sec:perf-measures} and Figure~\ref{fig:error-over-walltime} for further discussion.
\item[Reynolds Number $\boldsymbol{240}$.] Figure~\ref{fig:dl-3d-1-re240} and 
Table~\ref{tab:dl-3d-1-re240} show that DNN-MG captures the flow characteristics, as
demonstrated by enhanced resolution and the reproduction of features that were
unobservable at lower levels.
\item[Reynolds Number $\boldsymbol{300}$.] Figure~\ref{fig:dl-3d-1-re300} and 
Table~\ref{tab:dl-3d-1-re300} demonstrate that DNN-MG retains its effectiveness for Reynold number 300, with drag and lift predictions often closest to
reference solutions. Comparisons across
Tables~\ref{tab:dl-3d-1-re180},~\ref{tab:dl-3d-1-re240},
and~\ref{tab:dl-3d-1-re300} demonstrate the robust performance of DNN-MG, which
often outperform MG(4) solutions in terms of drag and lift accuracy.
\end{description}
In Figure~\ref{fig:vp-errors-1obs} we present the relative velocity and
pressure errors
\begin{equation}\label{eq:relerrors}
e_{\vec v} = \frac{\big\Vert (\tilde{\vec v}_n^{L+2} + d_n^{L+2} ) - \vec v_n^{L+2} \big\Vert_2}{\big\Vert \vec v_n^{L+2}  \Vert_2},\qquad
e_p = \frac{\big\Vert \tilde p_n^{L+2}-p_n^{L+2} \big\Vert_2}{\big\Vert p_n^{ L+2} \big\Vert_2},
\end{equation}
where $\Vert \cdot \Vert_2$ denotes the Euclidean norm of the discrete
solution.
The results show good improvement in terms of the accuracy of
DNN-MG over the coarse mesh solutions. Further, DNN-MG even outperforms the
solution at the intermediate level 4 between the coarse and the
reference level. In Figure~\ref{fig:lspec1obs} we plot the spectrum of the lift
for different Reynolds numbers. Again, DNN-MG is at least on par with the intermediate
level 4.

\subsection{Performance measurements}
\label{sec:perf-measures}
\begin{figure}
  \includegraphics{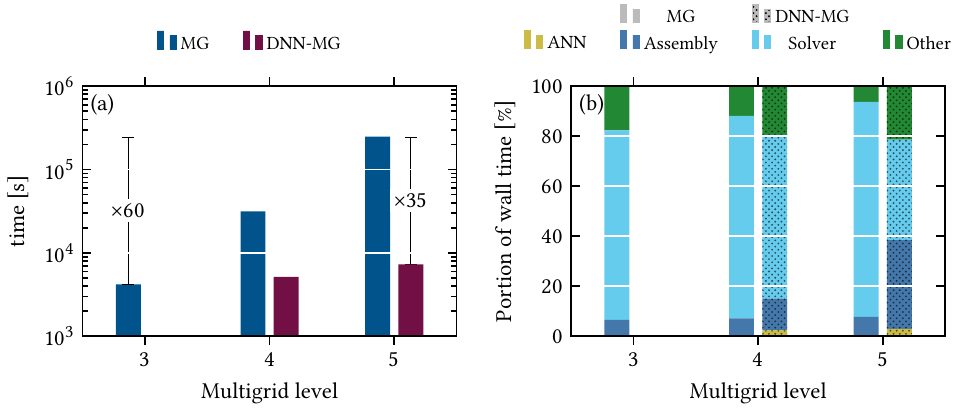}
  \caption{On the left the runtime of DNN-MG and MG is plotted over multiple
    levels of refinement. We see that DNN-MG has less than $50\%$ overhead
    compared to a coarse mesh simulation. On the right we see the contributions
    of different parts of the program to the wall clock time in percent. In DNN-MG($3+J$), the time spent in the solver (indicated in cyan) equals the time spent in the solver in MG($3$). We
    observe that the evaluation of the ANN takes less than $3\%$ of the
    time, which includes pre- and postprocessing of the in- and
    output.}\label{fig:perf}
  \vspace{.5em}
  \includegraphics{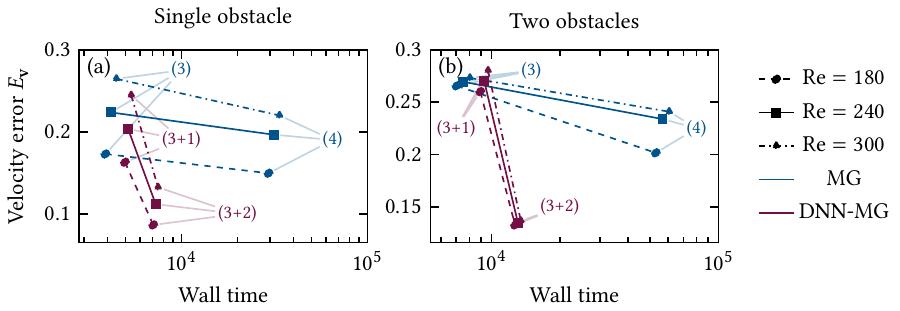}
  \caption{We plot the relative velocity error integrated over time
    \(\displaystyle E_{\vec v}=\paran{\frac{
      \int_0^8\lVert (\tilde{\vec v}_n^{L+2} + d_n^{L+2} ) - \vec v_n^{L+2}
      \rVert_2^2 \drv t}
          {\int_0^8\lVert \vec v_n^{L+2}  \rVert_2^2\drv t}}^{\frac{1}{2}}\)
    over the wall time of DNN-MG and MG.
    On the left, we show the result for the channel with one obstacle and on the
    right we show the result with two obstacles. We clearly see the potential of
    DNN-MG for two level predictions. One level predictions on the other hand
    barely improve the solution.}\label{fig:error-over-walltime}
\end{figure}
For the DNN-MG method to be advantageous, it must achieve two goals: Firstly, it
must enhance the accuracy compared to the solution on a coarse mesh. Secondly,
it must reduce the computation time compared to a solution on a fine mesh.
However, the Newton-Krylov geometric multigrid method already is of optimal
complexity~\cite{kimmritzParallelMultigrid2011,failerNewtonMultigridFramework2020}.
We note that for the simulation of incompressible flows, there exist one level
methods that can outperform multigrid methods in terms of efficiency, as
highlighted in~\autocite{ahmedAssessmentSolversSaddle2018}. Therefore, there may
be more efficient solvers than the Newton-Krylov geometric multigrid method in practice.
However, DNN-MG could still be applied to these one level methods, as we only
need to be able to refine the coarse grid.

In~\cite{margenbergNeuralNetworkMultigrid2022} the demonstrated efficiency of the numerical solver raised
the question what advantage DNN-MG offers, and we addressed it there for the 2D
case. The results indicated significant accuracy improvements and error
reductions in drag and lift calculations, and modest increases of the
wall time over coarser solutions. Notably, DNN-MG demonstrated substantial
generalization capabilities, achieving improved solutions for variants of channel
flows with different obstacle scenarios, including one with no obstacle at all.
Limitations were observed in the L-shaped domain, suggesting the need for more
sophisticated neural networks and diverse training data.

In~\cite{margenbergDeepNeuralNetworks2021}, we observed in a scalability analysis w.\,r.\,t.\
the neural network size that larger networks improve the predictions.
Specifically, larger models accurately replicated flow
frequencies in high-fidelity simulation and achieved good velocity and pressure error reductions with minimal
additional runtime costs. We also observed that DNN-MG improves the initial
guesses in Newton's method, which even reduced the runtime below the coarse mesh
solution.

Here, we adapt our theoretical findings obtained in the 2D case to the 3D one.
For this investigation, we rely on the properties of the Newton-Krylov geometric
multigrid method. If alternative solvers were used, the improvements in
efficiency (both theoretical and practical) would be different. Subsequently, we
analyze empirically whether the efficiency and scalability achieved in 2D translates to the
3D situation.

\paragraph{Complexity analysis}
The Newton-Krylov
geometric multigrid method has linear complexity \(\mathcal{O}(N)\) in the
number of DoFs \(N\). For each global mesh refinement, the number of DoFs
increase by a factor of 8, i.\,e.\ \(N^{L+1}\approx 8N^L\). The constant hidden in
\(\mathcal{O}(N)\), however, can be very significant, since on average 5 Newton
steps are required in each time step and within each Newton step one has to
compute on average 12 GMRES steps with one sweep of the geometric multigrid
solver as preconditioning. In addition, we perform 6 pre-smoothing and 6 post-smoothing steps in the geometric multigrid method. Thus, on each mesh
level \(L,\,L-1,\,\dots,\,0\), a total of approximately 700 Vanka smoothing
steps must be performed. One Vanka step thereby requires the inversion of about
\(\mathcal{O}(N^L/8)\) small matrices of size \(108\times 108\), one on each
element of the mesh, resulting in approximately \(108^3\approx \num{1259712}\)
operations. 
Since the complexity across all mesh levels sums up to about
\(N^L + N^{L-1}+\dots N_0 \approx N^L(1+8^{-1}+\cdots + 8^{-L})\approx
\frac{8}{7}N^L\) we can estimate the effort for the complete solution process on
level \(L\) as
\(5\cdot 12\cdot (6+6)\cdot 108^3\cdot \frac{8}{7}\approx 10^9 N^L\). We thereby
only counted the effort for smoothing and neglected all further matrix, vector
and scalar products, since these would increase the costs further and hence shift the results
of the analysis even more in favor of DNN-MG. If we were to solve the problem on mesh level \(L+1\), the
required effort would increase by a factor of 8 to approximately
\(\approx 8\cdot 10^9 N^L\) operations.

The DNN-MG method only requires the prolongation of the
solution to the next finer mesh and the evaluation of the neural network on each
patch. If we again consider patches of minimal size (one patch is one element of
mesh layer \(L+1\)) about \(\mathcal{O}(N^{L+1}/8)=\mathcal{O}(N^L)\) patches
must be processed. The effort for the evaluation of the network can be estimated
by the number of trainable parameters with \(N_c\) inputs. The DNN-MG approach
requires only one evaluation of the network in each time step. The number of
trainable parameters is for all network models of the order of single digit
millions. Hence, the effort for correcting the level \(L\) Newton-Krylov
solution by the neural network on level \(L+1\) can be estimated by
\(10^{7} N^L\), which is almost a factor of thousand lower than number of FLOPS
needed to solve on \(L+1\) (\(\approx 8\cdot10^9 N^L\)).

For the prediction of two levels, the estimate for the computational cost of
DNN-MG remains largely unchanged. While the size of the patches and the inputs
and outputs of the networks increase, the number of patches and network
evaluations do not change. 
Surprisingly, we also do not need to significantly increase
the network size and therefore give an upper bound of \(2\cdot 10^{7} N^L\)
trainable parameters. This is even more favorable than the prediction of a
single mesh level, as the numerical solution on level \(L+2\) incurs a much
higher cost of approximately \(64\cdot10^9 N^L\).

\paragraph{Computational performance of DNN-MG}
In Figure~\ref{fig:perf}\,(a) the runtimes of DNN-MG and a conventional multigrid
solution are plotted. For DNN-MG, the numerical solution level is constant with
$L=3$ and the number of predicted levels are $J=1,\:2$.
For the multigrid solution, the wall clock time increases by a factor of 8 with each added mesh level while for DNN-MG it increases at most by a factor of \(1.5\). For $\text{DNN-MG}(3+2)$ and
$\text{MG}(5)$ the wall clock time therefore differs by a factor of \(35\). 
In Figure~\ref{fig:error-over-walltime} the error accuracy is plotted over the
wall clock time. $\text{DNN-MG}(3+2)$ clearly has an advantage in terms of the
accuracy per runtime which once more underlines the computational efficiency of
DNN-MG.\@ Again, $\text{DNN-MG}(3+2)$ consistently outperforms MG(4). The advantage of $\text{DNN-MG}(3+1)$ is not as evident, which supports our statement that 2-level predictions are essential in 3D.

As can be seen in Figure~\ref{fig:perf}\,(b), the time spent in the solver in
DNN-MG is constant, which reduces the share of the solver in the wall time while
increasing the relative contribution of other parts of the program. When increasing $J$, i.e. the number of
predicted levels, the main increase in computation time is
due to the increased cost of the assembly of the right hand side and the
evaluation of the functionals, since these calculations are done on the finer
level. Note that these operations are also of order \(O(N)\), but with a much
smaller constant than other operations in the program, e.\,g.\ the operations
mentioned above in our theoretical considerations. When comparing MG($4$) and
DNN-MG($3+1$), or MG($5$) and DNN-MG($3+2$), the assembly cost
of DNN-MG is significantly lower. DNN-MG necessitates the assembly of a
residual and right-hand-side on level $L+J$ once per time step. Excluding
these two operations, the assembly times of DNN-MG($3+J$) are in the
same range as MG($3$). This makes the DNN-MG method
highly efficient for 2 level predictions (see also Figure~\ref{fig:error-over-walltime}). 

An investigation with more levels was not feasible in our work
since the version of Gascoigne3D we are using only support shared-memory parallelism
and we are therefore restricted in terms of the problem size. 
Lowering the base level led to too much inaccuracies for DNN-MG to be applicable. 
Further investigations of the scalability of the method when predicting more mesh-levels
is part of future work. Predicting weights of higher order polynomials instead
of the weights of basis functions on refined meshes could further reduce the
computational burden and increase the computational efficiency.

The evaluation of the DNN takes less than $3\%$ of the wall time including pre- and post-processing of the in- and output.
For a fixed $L$, this ratio is independent of the number of predicted levels $J$.
Considering that the networks are of similar size for $L+1$
and $L+2$ (cf.~Table~\ref{tab:orgef2be82}), this is expected.
This is in support of our previous remark that the
increase in runtime from DNN-MG($L+1$) to DNN-MG($L+2$) results from
the computational cost of assembling the residual and right-hand side on level $L+1$ and $L + 2$, respectively.

If $L$ is not fixed, the amount of data, and thus the cost of evaluating the neural network, increases by
a factor of 8 in 3D (4 in 2D), resulting in identical scaling as the numerical solver. Thus,
the ratio of 3\% will not change significantly for different $L$. We note that this is less
than our previous results~\cite{margenbergNeuralNetworkMultigrid2022} in 2D.
This is due to substantial improvements of our implementation. In particular, we now collect the
input of the neural network and process the output is parallel and the
evaluation of the network is done on a GPU.\@  For the parallelization, we use
\texttt{OpenMP} directives. For evaluation of the preprocessing operator
$\mathcal{Q}$ (cf.~\eqref{eq:defQ}), we utilize a simple parallel for loop over the patches. For evaluation of the
post-processing operator $\mathcal{L}$ (cf.~\eqref{eq:defL}), we again use a
parallel for loop over the patches and protect write access to the global vector
by a \texttt{\#pragma omp atomic update}. Through the combination of CPU and GPU computing, DNN-MG not only reduces
the cost of accurate simulations, but it is also well-suited for heterogeneous
computing environments.

\subsection{Generalizability}
\label{sec:generalize}
\begin{figure}[t]
  \includegraphics{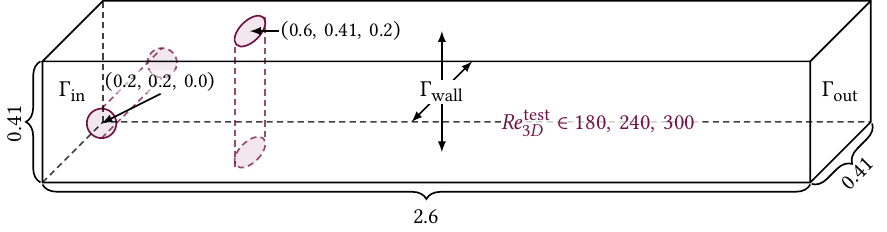}
  \caption{Geometry of the training and test scenario with a parabolic inflow
    profile at $\Gamma_{\textrm{in}}$, do-nothing boundary conditions at the
    outflow boundary $\Gamma_{\textrm{out}}$ and no-slip conditions on the walls
    $\Gamma_{\textrm{wall}}$. The center of the obstacle is at $(0.5,\,0.2)$.}\label{fig:2obs-geo}
  \end{figure}

In order to ensure the practicality of the DNN-MG method, it must generalize
well to similar flows beyond those seen during training. Previous results on the
3D benchmark demonstrate the network's ability to do so under small
geometric perturbations and varying Reynolds number. Here, we evaluate the
network's performance under more substantial changes considering a channel with two
obstacles shown in Figure~\ref{fig:2obs-geo}. As for the single obstacle case,
we test the performance under varying Reynolds numbers with a viscosity of
\(\nu=\num{5.e-4}\) (\(Re=240\)), \(\nu=\num{4e-4}\) (\(Re=300\)) and
\(\nu=\num{6.67e-4}\) (\(Re=300\)). We adopt some settings from the single
obstacle case. For instance, the boundary conditions, in particular the inflow
at the left boundary \(\Gamma_{\text{in}}\) given by~\eqref{eq:inflow}, the
initial velocity \(\vec{v}(0) = \vec{0}\), and we chose the time step size
\(k=\num{0.008}\) on the interval \(I=(0,8]\).

\begin{table}[t]
  \caption{Spatial meshes with increasing number of
    multigrid levels and unknowns.}
    \label{tab:org6b4401b}
  \centering
  \begin{tabular}{llrr}
    \toprule
                              & Level & \# DoF             & Patch size \\[0pt]
    \midrule
                       Coarse & MG(\(3\)) & \(\num{533148}\)   & 0          \\[0pt]
    Reference for $\text{DNN-MG}(3+1)$ & MG(\(4\)) & \(\num{4097340}\)  & 8         \\[0pt]
    Reference for $\text{DNN-MG}(3+2)$ & MG(\(5\)) & \(\num{32115324}\) & 64        \\[0pt]
    \bottomrule
  \end{tabular}
\end{table}

\begin{figure}[htbp]
  \centering%
  \includegraphics{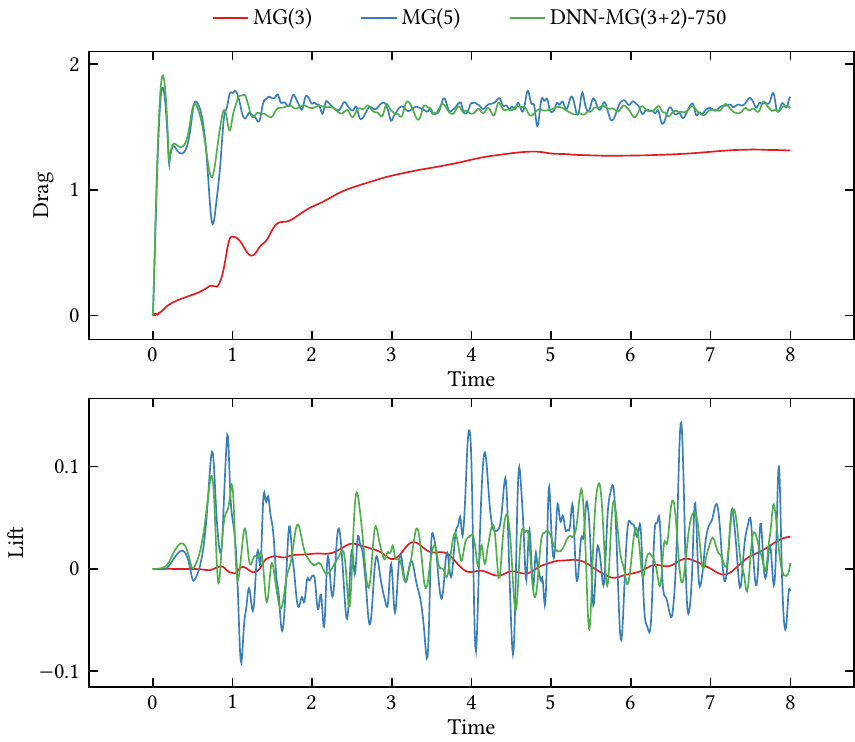}
  \caption{The drag (left) and lift (right) functionals for the channel with two
    obstacles at $\text{Re}=180$ are shown. The results are presented for three
    discretizations: The coarse mesh $\text{MG}(3)$, the reference
    $\text{MG}(5)$, and the $\text{DNN-MG}(3+2)$ improved by a neural
    network.}\label{fig:dl-3d-2-re180}
  \captionof{table}{The $\min$, $\max$, mean and amplitude of drag and lift
    functionals in the same setting as above in Figure~\ref{fig:dl-3d-2-re180}.
    The best results in comparison to the $\text{MG}(5)$ solution are
    highlighted with color coding: \bestf{Reference}, \bests{1st best},
    \bestt{2nd best}.}\label{tab:dl-3d-2-re180}
  \begin{tabular}{lrrrrrrrr}
    \toprule
    Type & \(\min C_{d}\) & \(\max C_{d}\) & \(\widebar{C_{d}}\) & amp \(C_{d}\) & \(\min C_{l}\) & \(\max C_{l}\) & \(\widebar{C_{l}}\) & amp \(C_{l}\)\\[0pt]
    \midrule
    $\text{MG}(3)$ & 0.7674 & 1.3213 & 1.2176 & 0.5539 & -0.0087 & 0.0313 & \bests{0.0097} & 0.0400\\[0pt]
    $\text{MG}(4)$ & 0.7905 & 1.5601 & 0.7696 & 1.3747 & -0.0091 & 0.0346 & 0.0061 & 0.0404\\[0pt]
    $\text{MG}(5)$ & \bestf{1.5089} & \bestf{1.7903} & \bestf{1.6636} & \bestf{0.2814} & \bestf{-0.0872} & \bestf{0.1428} & \bestf{0.0097} & \bestf{0.2299}\\[0pt]
    $\text{DNN-MG}(3+1)$-512 & \bests{1.5179} & \bests{1.7259} & \bests{1.6342} & \bests{0.2079} & -0.0487 & \bests{0.0907} & 0.0280 & \bestt{0.1394}\\[0pt]
    $\text{DNN-MG}(3+2)$-512 & 1.6038 & 1.6457 & 1.6265 & 0.0419 & -0.0432 & -0.0112 & -0.0237 & 0.0320\\[0pt]
    $\text{DNN-MG}(3+1)$-750 & 1.5980 & 1.6897 & 1.6436 & 0.0917 & \bestt{-0.0545} & 0.0784 & 0.0183 & 0.1330\\[0pt]
    $\text{DNN-MG}(3+2)$-750 & \bestt{1.5591} & \bestt{1.7051} & \bestt{1.6376} & \bestt{0.1460} & \bests{-0.0591} & \bestt{0.0834} & \bestt{0.0096} & \bests{0.1426}\\[0pt]
    \bottomrule
  \end{tabular}
\end{figure}

\begin{figure}[htbp]
  \centering%
  \includegraphics{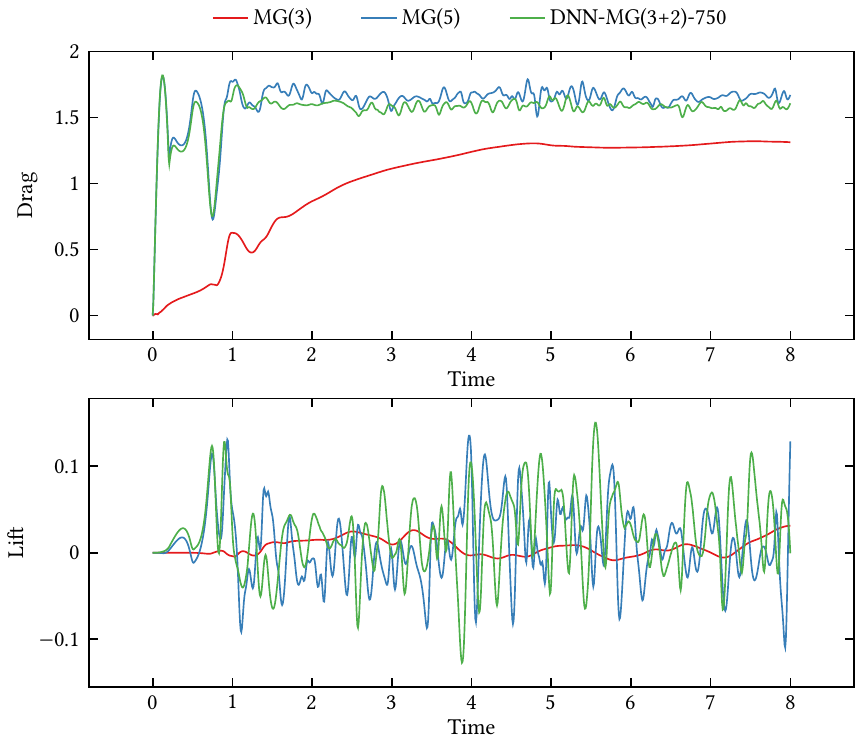}
  \caption{The drag (left) and lift (right) functionals for the channel with two
    obstacles at $\text{Re}=240$ are shown. The results are presented for three
    discretizations: The coarse mesh $\text{MG}(3)$, the reference
    $\text{MG}(5)$, and the $\text{DNN-MG}(3+2)$ improved by a neural
    network.}\label{fig:dl-3d-2-re240}
  \captionof{table}{The $\min$, $\max$, mean and amplitude of drag and lift
    functionals in the same setting as above in Figure~\ref{fig:dl-3d-2-re240}.
    The best results in comparison to the $\text{MG}(5)$ solution are
    highlighted with color coding: \bestf{Reference}, \bests{1st best},
    \bestt{2nd best}.}\label{tab:dl-3d-2-re240}
  \begin{tabular}{lrrrrrrrr}
    \toprule
    Type & \(\min C_{d}\) & \(\max C_{d}\) & \(\widebar{C_{d}}\) & amp \(C_{d}\) & \(\min C_{l}\) & \(\max C_{l}\) & \(\widebar{C_{l}}\) & amp \(C_{l}\)\\[0pt]
    \midrule
    $\text{MG}(3)$ & 0.7674 & 1.3213 & 1.2176 & 0.5539 & -0.0087 & 0.0313 & \bests{0.0097} & 0.0400\\[0pt]
    $\text{MG}(4)$ & 0.7905 & 1.5601 & 0.7696 & 1.3747 & -0.0091 & 0.0346 & 0.0061 & 0.0404\\[0pt]
    $\text{MG}(5)$ & \bestf{1.5090} & \bestf{1.7907} & \bestf{1.6642} & \bestf{0.2818} & \bestf{-0.1577} & \bestf{0.1693} & \bestf{0.0083} & \bestf{0.3270}\\[0pt]
    $\text{DNN-MG}(3+1)$-512 & 1.3311 & 1.6166 & 1.3971 & \bests{0.2855} & -0.0505 & 0.0960 & 0.0291 & 0.1466\\[0pt]
    $\text{DNN-MG}(3+2)$-512 & 1.4353 & 1.5870 & 1.5126 & 0.1517 & \bests{-0.1398} & 0.0790 & -0.0528 & 0.2188\\[0pt]
    $\text{DNN-MG}(3+1)$-750 & \bestt{1.4527} & \bestt{1.6565} & \bestt{1.5547} & 0.2038 & -0.0832 & \bestt{0.1507} & \bestt{0.0219} & \bestt{0.2339}\\[0pt]
    $\text{DNN-MG}(3+2)$-750 & \bests{1.4644} & \bests{1.6823} & \bests{1.5929} & \bestt{0.2179} & \bestt{-0.1325} & \bests{0.1684} & 0.0279 & \bests{0.3009}\\[0pt]
    \bottomrule
  \end{tabular}
\end{figure}

\begin{figure}[htbp]
  \centering%
  \includegraphics{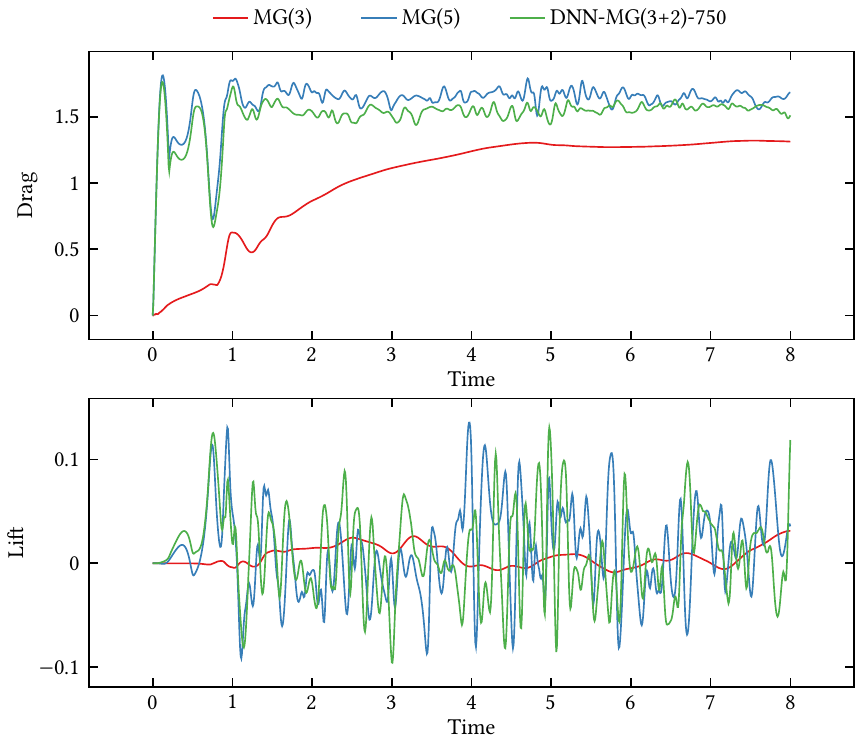}
  \caption{The drag (left) and lift (right) functionals for the channel with two
    obstacles at $\text{Re}=300$ are shown. The results are presented for three
    discretizations: The coarse mesh $\text{MG}(3)$, the reference
    $\text{MG}(5)$, and the $\text{DNN-MG}(3+2)$ improved by a neural
    network.}\label{fig:dl-3d-2-re300}
  \captionof{table}{The $\min$, $\max$, mean and amplitude of drag and lift
    functionals in the same setting as above in Figure~\ref{fig:dl-3d-2-re300}.
    The best results in comparison to the $\text{MG}(5)$ solution are
    highlighted with color coding: \bestf{Reference}, \bests{1st best},
    \bestt{2nd best}.}\label{tab:dl-3d-2-re300}
  \begin{tabular}{lrrrrrrrr}
    \toprule
    Type & \(\min C_{d}\) & \(\max C_{d}\) & \(\widebar{C_{d}}\) & amp \(C_{d}\) & \(\min C_{l}\) & \(\max C_{l}\) & \(\widebar{C_{l}}\) & amp \(C_{l}\)\\[0pt]
    \midrule
    $\text{MG}(3)$ & 0.7674 & 1.3213 & 1.2176 & 0.5539 & -0.0087 & 0.0313 & \bests{0.0097} & 0.0400\\[0pt]
    $\text{MG}(4)$ & 0.7905 & 1.5601 & 0.7696 & 1.3747 & -0.0091 & 0.0346 & 0.0061 & 0.0404\\[0pt]
    $\text{MG}(5)$ & \bestf{1.5092} & \bestf{1.7908} & \bestf{1.6568} & \bestf{0.2816} & \bestf{-0.0871} & \bestf{0.1354} & \bestf{0.0122} & \bestf{0.2225}\\[0pt]
    $\text{DNN-MG}(3+1)$-512 & 1.3049 & 1.5990 & 1.5111 & \bests{0.2941} & -0.0471 & \bestt{0.1589} & 0.0317 & 0.2060\\[0pt]
    $\text{DNN-MG}(3+2)$-512 & 1.3021 & 1.5548 & 1.3596 & \bestt{0.2526} & -0.1398 & 0.0790 & -0.0528 & \bestt{0.2188}\\[0pt]
    $\text{DNN-MG}(3+1)$-750 & \bestt{1.4050} & \bestt{1.6214} & \bestt{1.4465} & 0.2164 & \bestt{-0.0684} & 0.2384 & 0.1214 & 0.3068\\[0pt]
    $\text{DNN-MG}(3+2)$-750 & \bests{1.4387} & \bests{1.6319} & \bests{1.5499} & 0.1932 & \bests{-0.0957} & \bests{0.1343} & \bestt{0.0067} & \bests{0.2300}\\[0pt]
    \bottomrule
  \end{tabular}
\end{figure}

Since we are interested in understanding the generalization abilities of DNN-MG, we use the network trained for the single-obstacle channel at
Reynolds number \(Re=200\). The neural network
component of DNN-MG operates on level \(5\) as fine resolution and level \(3\)
as coarse resolution in this section. In order to not overburden this paper we
leave out the detailed results of a single level prediction. The complexity of the finite element meshes is given in Table~\ref{tab:org6b4401b}.
\begin{description}[style=unboxed,leftmargin=0cm]\itemsep1pt \parskip0pt \parsep0pt
\item[Channel with two obstacles at $\text{Re}=\boldsymbol{180}$.]
DNN-MG shows good prediction accuracy in flow dynamic and
significantly improves drag and lift coefficient accuracy and closely
reconstructs the dynamics of the reference solution, despite temporal
discrepancies, see Figure~\ref{fig:dl-3d-2-re180} and Table~\ref{tab:dl-3d-2-re180}.
As mentioned above, it is important to use DNN-MG with two-level predictions,
because single level predictions lack reliability. This is even more pronounced
for the two obstacle generalization test.
\item[Channel with two obstacles at $\text{Re}=\boldsymbol{240}$.] 
DNN-MG continues to improve the accuracy and effectively
captures untrained features of the solution, see Figure~\ref{fig:dl-3d-2-re240} and 
Table~\ref{tab:dl-3d-2-re240}. Figure~\ref{fig:fields-2obs-re240} shows the
capability of the method to increase the resolution of the flow and
to add features that could not be observed on lower levels but which are present on
higher levels. This emphasizes the power of DNN-MG's localized predictions.
\item[Channel with two obstacles at $\text{Re}=\boldsymbol{300}$.]
In Figure~\ref{fig:dl-3d-2-re300} and
Table~\ref{tab:dl-3d-2-re300} we show the drag and lift forces predicted by
DNN-MG.\@ The results underscore DNN-MG's generalization capabilities. They are
consistent across different Reynolds numbers, with results at $\text{Re}=300$
analogous to those at $\text{Re}=240$ and $\text{Re}=180$.
It should be emphasized in particular that the second obstacle is hit by an already fully developed non-stationary flow. This situation does not exist at all in the one-obstacle case and thus in the training data.
From the above observations we conclude that
DNN-MG offers great potential in terms of computational efficiency.
\end{description}

\begin{figure}
  \centering \includegraphics{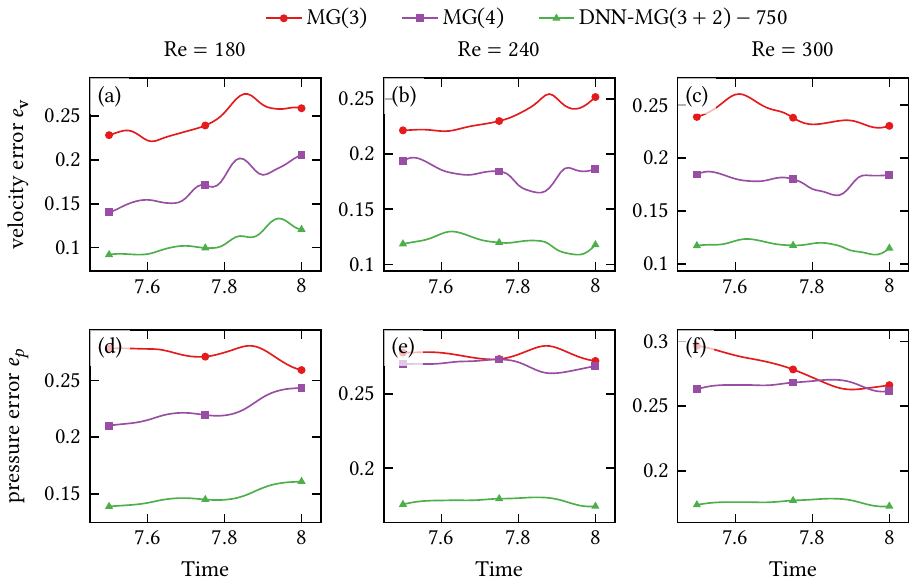}
  \caption{Relative velocity errors (first row) and pressure errors (second row)
    for the channel with two
    obstacles. We compare coarse solutions with and without ANN correction ($\text{MG}(3)$ and
    $\text{DNN-MG}(3+2)-750$) to the reference solution on level $L+2$
    for different Reynolds numbers. As in the case with a single obstacle, we
    show the errors of the $\text{MG}(4)$ solution as well, which we again
    outperform.}\label{fig:vp-errors-2obs}
  \vspace{.5em}%
  \includegraphics{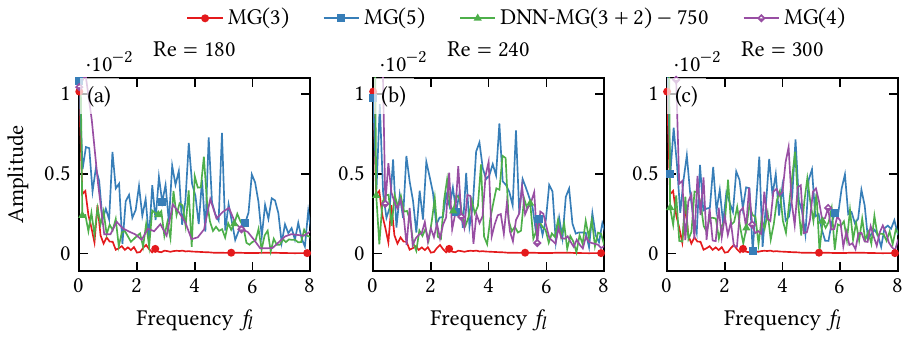}%
  \caption{The spectrum of the lift functional for the channel flow with two
    obstacles at different Reynolds numbers. We plot the results of the coarse
    mesh $\text{MG}(3)$, the reference $\text{MG}(5)$, the intermediate level
    $\text{MG}(4)$ and the $\text{DNN-MG}(3+2)-750$ improved by a neural
    network.}\label{fig:lspec2obs}
\end{figure}

Figure~\ref{fig:vp-errors-2obs} demonstrates the good performance of DNN-MG
in terms of velocity and pressure error. The figure
shows a significant improvement of the accuracy over an $\text{MG}(3)$ solution and
remarkably even outperforms the $\text{MG}(4)$ solution. This further
substantiates the efficiency in terms of runtime that we demonstrated in
Section~\ref{sec:perf-measures}. In Figure~\ref{fig:lspec2obs} we plot the spectrum of the lift
for different Reynolds numbers. Again, DNN-MG is at least on par with the intermediate
level 4.

Overall, DNN-MG is robust with respect to changes of the geometry and material
parameters and is able to improve solutions to PDEs in terms of characteristic
quantities of interest and error measures.

\section{Conclusion}
\label{sec:conc}
We have presented the deep neural network multigrid solver (DNN-MG) which uses a
deep neural network to improve the efficiency of a classical finite element solver,
e.\,g.\ for the simulation of the Navier-Stokes equations in 3D.
The grids we use are unstructured and pre-adapted to the problem.
Previously demonstrated
for 2D simulations, we extended DNN-MG to 3D and reformulated it
in a rigorous manner.
Despite the increased complexity of 3D flows, the algorithm remained applicable and
delivered even larger speed-ups compared to 2D while retaining its efficiency,
generalizability, and scalability.

We established the efficacy of DNN-MG for large-scale simulations in regimes
where direct solvers are not feasible anymore. The overhead is small and in 3D
we can accelerate high fidelity simulations by a factor of 35, although with a trade-off concerning the quality of the solution.
However, using predictions across two levels, we
were able to consistently outperform classical solutions on the intermediate
level both in terms of solution quality and wall time.
This efficiency-performance trade-off establishes DNN-MG as a highly
promising approach for accelerating numerical solution methods for PDEs.

We also demonstrated the generalization capabilities of DNN-MG which stem largely from its local
approach. We trained the neural network only with a single flow scenario at
$\text{Re}=200$ and got significantly improved results for lower and higher Reynolds numbers. 
In particular, our results showed that DNN-MG substantially improves the solution accuracy, as
measured by the $l^{2}$-error, and reduces the errors of the drag and lift
functionals across all tested scenarios.
The results for two-obstacle flow scenarios that DNN-MG is able to successfully
improve the solution in situations that are not present in the training data
across a range of Reynolds numbers.

Especially noteworthy is that DNN-MG is able to predict flow profiles of completely different dynamics. The coarse grid MG(\(3\)) and the fine MG(\(5\)) grid considered by us produce flows with very different character. Thus, especially in the test problem with two obstacles, there are hardly any oscillations on the coarse grid. These are, nonetheless, correctly reproduced by DNN-MG. The local approach is thus able to identify global structures of the Navier-Stokes solution and to correct the solution in both local and temporal dynamics.

Although DNN-MG generalizes well, there are limits to this approach. In future
work, we want to develop an online learning approach to retrain the network
adaptively based on the uncertainty of the predictions. Further, we want to
incorporate physical information into the neural network by including the residual of the PDE into the
loss function. We also want to investigate the application of DNN-MG to other
PDEs.

\FloatBarrier%

\subsection*{Acknowledgement}
NM acknowledges support by the Helmholtz-Gesellschaft grant number HIDSS-0002
DASHH. Computational resources (HSUper) were provided by the project
hpc.bw, funded by dtec.bw — Digitalization and Technology Research Center of the
Bundeswehr. dtec.bw is funded by the European Union – NextGenerationEU.
The work of TR was supported by the Deutsche Forschungsgemeinschaft (DFG, German Research Foundation) - 314838170, GRK 2297 MathCoRe.

\printbibliography
\appendix
\section{Descriptive Statistics of the data from the test cases}\label{app:stat}
\begin{minipage}{\textwidth}
  \captionof{table}{Descriptive statistics of the 6 test cases.}\label{tab:stats1}
  \small
  \sisetup{round-mode=places,round-precision=2,exponent-mode=scientific,print-zero-exponent=false,tight-spacing=true}
  \centering
  \begin{minipage}{0.49\textwidth}
    \begin{center}
      {\bfseries One Obstacle, Re=180}\\[.5em]
      \begin{tabular}{crrrr}
        \toprule
        & max & mean & min & $\sigma$ \\[0pt]
        \midrule
        $\vec v_{x}$ & \num{2.718} & \num{0.979}     & -\num{0.942}& \num{0.705} \\[0pt]
        $\vec v_{y}$ & \num{1.663} & \num{0.0005716}  & -\num{1.6788} & \num{0.2694} \\[0pt]
        $\vec v_{z}$ & \num{1.275} & \num{3.81206e-11}  & -\num{1.275} & \num{0.172} \\[0pt]
        $p$          & \num{3.398} & \num{0.24343}     & -\num{1.523} & \num{0.5424} \\[0pt]
        \bottomrule
      \end{tabular}
    \end{center}
  \end{minipage}
  \begin{minipage}{0.49\textwidth}
    \begin{center}
      {\bfseries Two Obstacles, Re=180}\\[.5em]
      \begin{tabular}{crrrr}
        \toprule
        & max & mean & min & $\sigma$ \\[0pt]
        \midrule
        $\vec v_{x}$ & \num{2.768} & \num{0.987}     & -\num{1.737}  & \num{0.7444} \\[0pt]
        $\vec v_{y}$ & \num{1.656} & \num{-0.0005781} & -\num{1.607}   & \num{0.248} \\[0pt]
        $\vec v_{z}$ & \num{1.928} & \num{4.912e-13} & -\num{2.049}  &  \num{0.2782} \\[0pt]
        $p$          & \num{4.05} & \num{0.2932}    & -\num{1.724}  & \num{0.5577} \\[0pt]
        \bottomrule
      \end{tabular}
    \end{center}
  \end{minipage}

  \vspace{.5em}
  \begin{minipage}{0.49\textwidth}
    \begin{center}
      {\bfseries One Obstacle, Re=240}\\[.5em]
      \begin{tabular}{crrrr}
        \toprule
        & max & mean & min & $\sigma$ \\[0pt]
        \midrule
        $\vec v_{x}$ & \num{2.823} & \num{0.9798}     & -\num{1.2}& \num{0.7059} \\[0pt]
        $\vec v_{y}$ & \num{1.729} & \num{0.003113}  & -\num{1.768} & \num{0.2851} \\[0pt]
        $\vec v_{z}$ & \num{1.6} & \num{2.433e-11}  & -\num{-1.6} & \num{0.1911} \\[0pt]
        $p$          & \num{3.265} & \num{0.1351}     & -\num{2.082} & \num{0.5277} \\[0pt]
        \bottomrule
      \end{tabular}
    \end{center}
  \end{minipage}
  \begin{minipage}{0.49\textwidth}
    \begin{center}
      {\bfseries Two Obstacles, Re=240}\\[.5em]
      \begin{tabular}{crrrr}
        \toprule
        & max & mean & min & $\sigma$ \\[0pt]
        \midrule
        $\vec v_{x}$ & \num{2.716} & \num{0.9852}     & -\num{-1.763}  & \num{0.7442} \\[0pt]
        $\vec v_{y}$ & \num{2.019} & \num{0.001095} & -\num{2.03}   & \num{0.2522} \\[0pt]
        $\vec v_{z}$ & \num{2.032} & \num{0.003468} & -\num{2.057}  &  \num{0.2813} \\[0pt]
        $p$          & \num{4.038} & \num{0.3004}    & -\num{1.949}  & \num{0.5466} \\[0pt]
        \bottomrule
      \end{tabular}
    \end{center}
  \end{minipage}

  \vspace{.5em}
  \begin{minipage}{0.49\textwidth}
    \begin{center}
      {\bfseries One Obstacle, Re=300}\\[.5em]
      \begin{tabular}{crrrr}
        \toprule
        & max & mean & min & $\sigma$ \\[0pt]
        \midrule
        $\vec v_{x}$ & \num{3.077} & \num{0.9817}     & -\num{1.422}& \num{0.6923} \\[0pt]
        $\vec v_{y}$ & \num{2.002} & \num{0.000295}  & -\num{1.777} & \num{0.2998} \\[0pt]
        $\vec v_{z}$ & \num{1.74} & \num{0.0000005667}  & -\num{1.74} & \num{0.2071} \\[0pt]
        $p$          & \num{3.149} & \num{0.07706}     & -\num{3.317} & \num{0.5265} \\[0pt]
        \bottomrule
      \end{tabular}
    \end{center}
  \end{minipage}
  \begin{minipage}{0.49\textwidth}
    \begin{center}
      {\bfseries Two Obstacles, Re=300}\\[.5em]
      \begin{tabular}{crrrr}
        \toprule
        & max & mean & min & $\sigma$ \\[0pt]
        \midrule
        $\vec v_{x}$ & \num{3.015} & \num{0.9679}     & -\num{1.784}  & \num{0.7415} \\[0pt]
        $\vec v_{y}$ & \num{2.121} & \num{-0.0004126} & -\num{2.06}   & \num{0.2683} \\[0pt]
        $\vec v_{z}$ & \num{2.079} & \num{-0.004395} & -\num{2.099}  &  \num{0.2889} \\[0pt]
        $p$          & \num{4.018} & \num{0.2809}    & -\num{1.621}  & \num{0.5804} \\[0pt]
        \bottomrule
      \end{tabular}
    \end{center}
  \end{minipage}
\end{minipage}
\section{Snapshots of the Flows}
\begin{minipage}{\textwidth}
  \centering
  \begin{tikzpicture}[text=white]
    \node[anchor=south west,inner sep=0] (A) at (0,0)
    {\includegraphics[width=\textwidth]{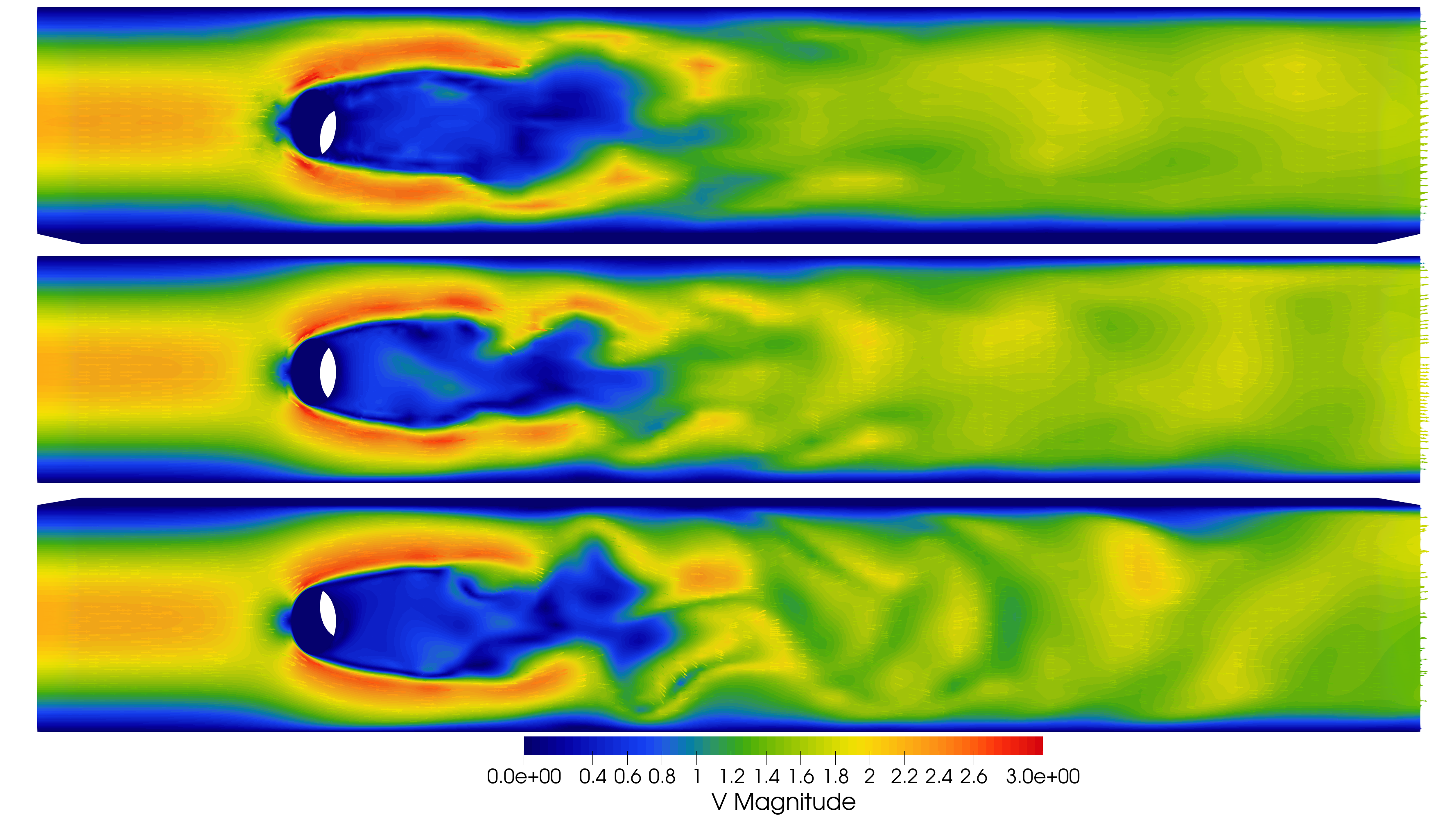}};
    \node[anchor=south west,inner sep=0] (Aa) at (0.5,8.6) {$\text{MG}(3)$};
    \node[anchor=south west,inner sep=0] (Ab) at (0.5,5.85) {DNN-MG($3+2$)};
    \node[anchor=south west,inner sep=0] (Ac) at (0.5,3.1) {$\text{MG}(5)$};
    \node[anchor=south west,inner sep=0] (B) [below=of A]
    {\includegraphics[width=\textwidth]{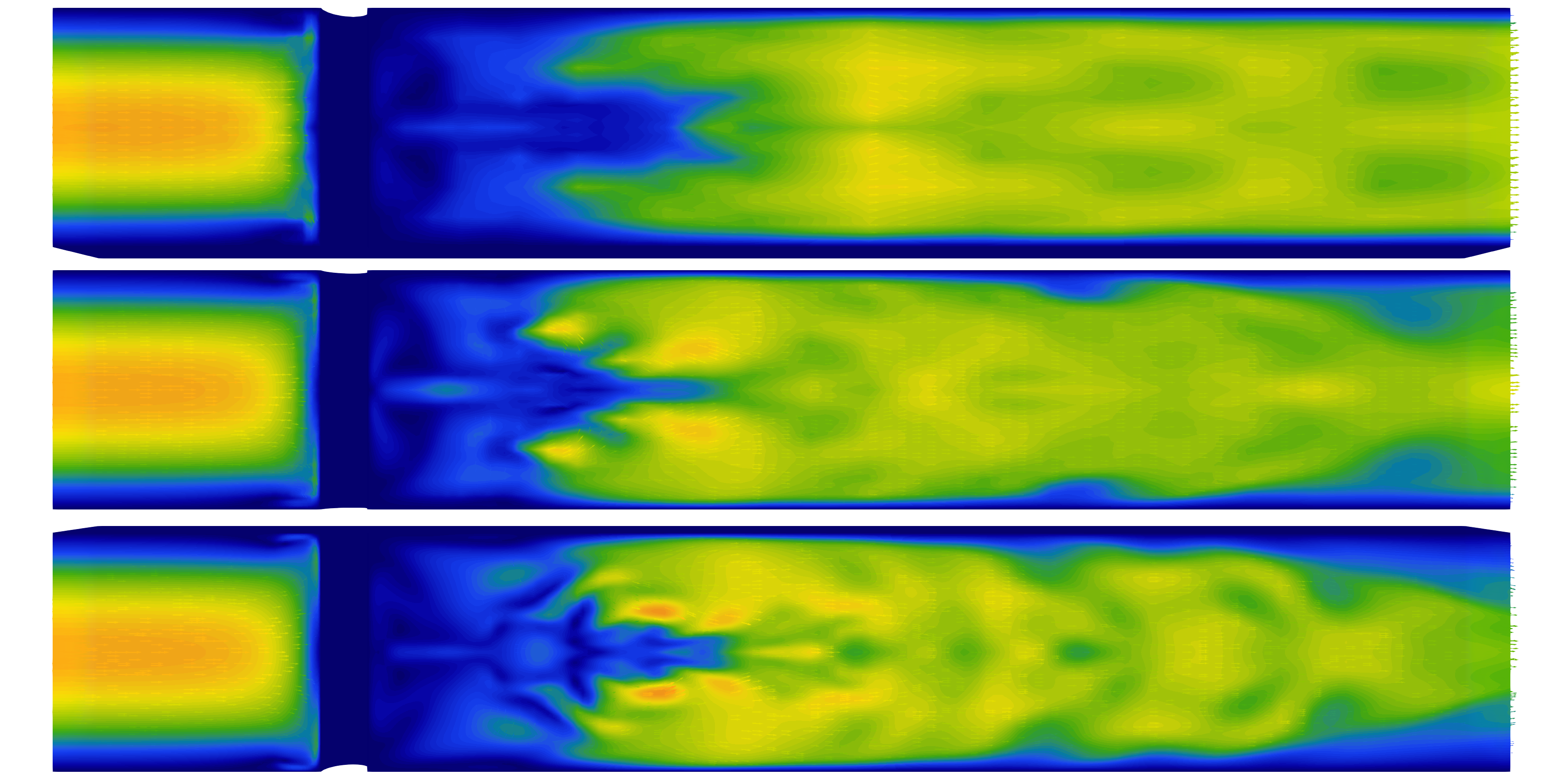}};
    \node[anchor=south west,inner sep=0] (Ba) at (0.6,-1.5) {$\text{MG}(3)$};
    \node[anchor=south west,inner sep=0] (Bb) at (0.6,-4.2)
    {$\text{DNN-MG}(3+2)$}; \node[anchor=south west,inner sep=0] (Bc) at
    (0.6,-6.9) {$\text{MG}(5)$};
  \end{tikzpicture}
  \captionof{figure}{Comparison of the magnitude of the velocity field for the channel
    with a cylindrical obstacle with an elliptical cross-section at
    $\mathrm{Re}=240$ and time $t=7.5$. The results are shown for
    $\text{DNN-MG}(3+2)$, along with the reference solution obtained with $5$
    levels ($\text{MG}(5)$), as well as the coarse mesh solution at level $3$
    ($\text{MG}(3)$). The $\text{DNN-MG}(3+2)$ method utilizes a network trained
    on simulation data from a round obstacle. The first three images show a
    cross-section of the domain in the $xy$-plane and the others show a
    cross-section of the domain in the $xz$-plane.}\label{fig:fields-1obs-re240}
\end{minipage}

\begin{minipage}{\textwidth}
  \centering
  \begin{tikzpicture}[text=white]
    \node[anchor=south west,inner sep=0] (A) at (0,0)
    {\includegraphics[width=\textwidth]{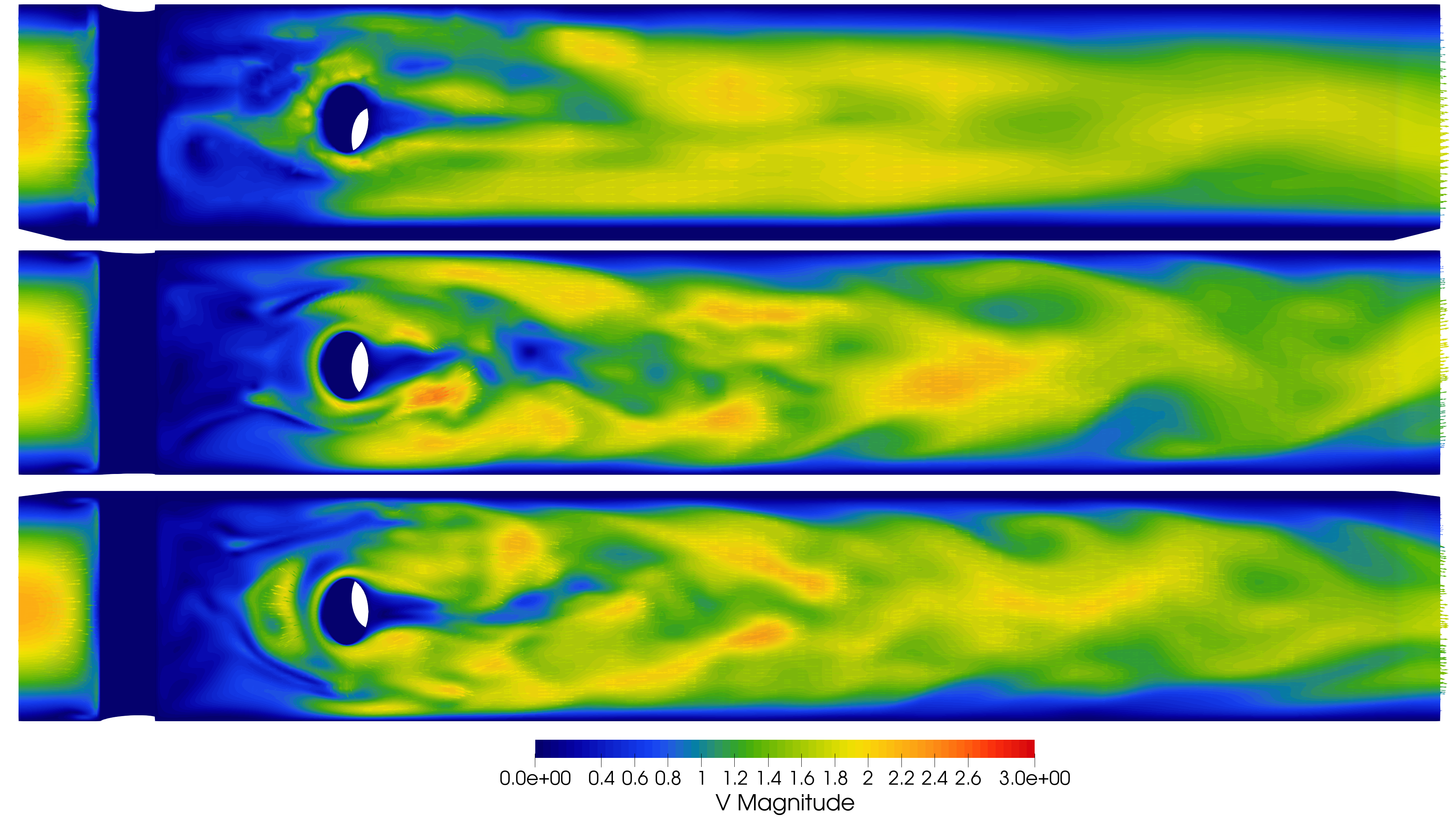}};
    \node[anchor=south west,inner sep=0] (Aa) at (0.3,8.65) {$\text{MG}(3)$};
    \node[anchor=south west,inner sep=0] (Ab) at (0.3,5.9)
    {$\text{DNN-MG}(3+2)$}; \node[anchor=south west,inner sep=0] (Ac) at
    (0.3,3.15) {$\text{MG}(5)$}; \node[anchor=south west,inner sep=0] (B)
    [below=of A] {
      \includegraphics[width=\textwidth]{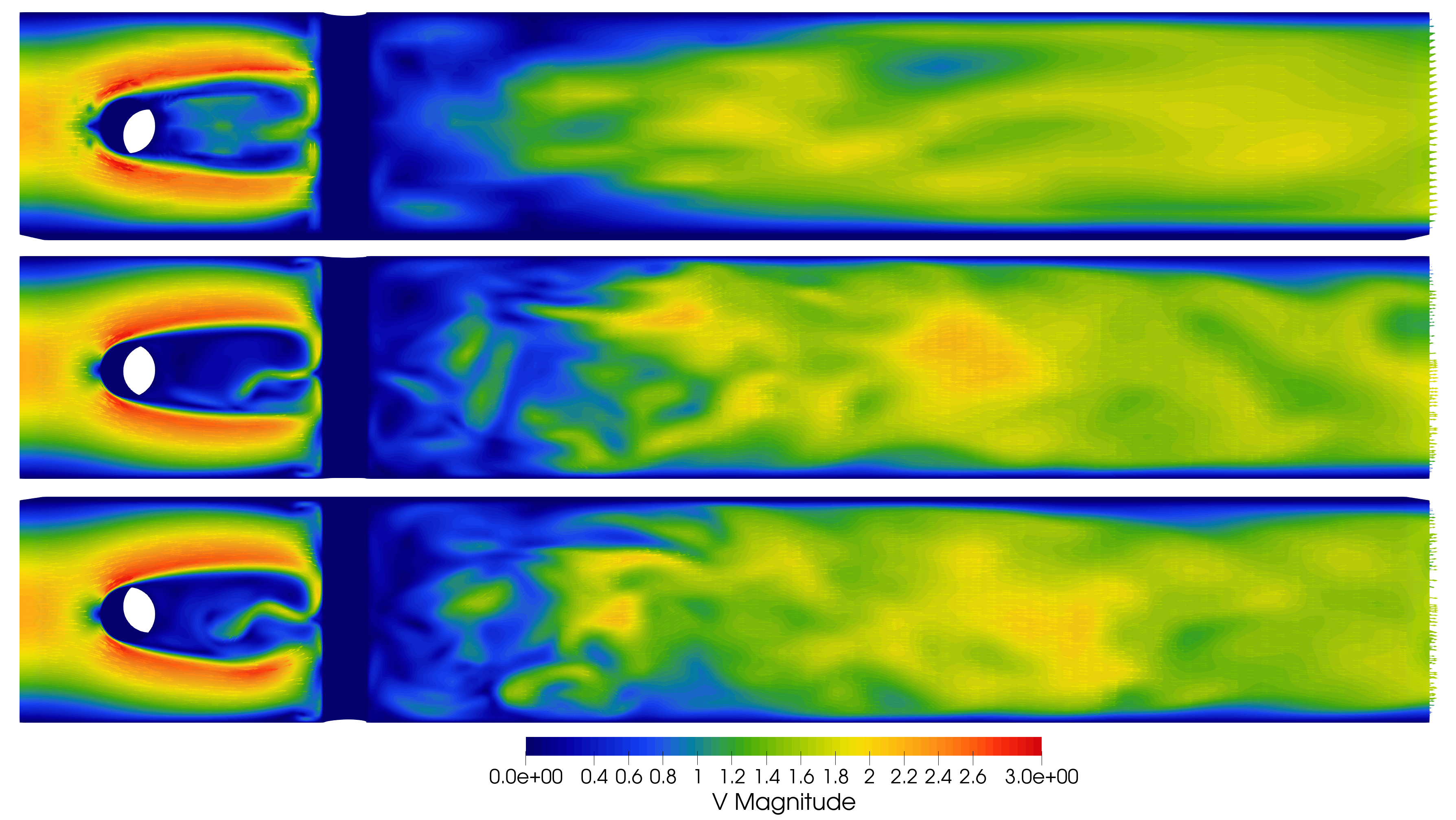}};
    \node[anchor=south west,inner sep=0] (Ba) at (0.3,-1.55) {$\text{MG}(3)$};
    \node[anchor=south west,inner sep=0] (Bb) at (0.3,-4.25)
    {$\text{DNN-MG}(3+2)$}; \node[anchor=south west,inner sep=0] (Bc) at
    (0.3,-6.95) {$\text{MG}(5)$};
  \end{tikzpicture}
  \captionof{figure}{Comparison of the magnitude of the velocity field for the channel
    with two cylindrical obstacles at time $t=7.5$. The first obstacle has a
    circular cross-section and the second obstacle has an elliptical
    cross-section, which results in a Reynolds number of $\mathrm{Re}=240$. The
    results are shown for $\text{DNN-MG}(3+2)$, along with the reference
    solution obtained with $5$ levels ($\text{MG}(5)$), as well as the coarse
    mesh solution at level $3$ ($\text{MG}(3)$). The $\text{DNN-MG}(3+2)$ method
    utilizes a network trained on simulation data from a round obstacle. The
    first three images show a cross-section of the domain in the $xy$-plane and
    the others show a cross-section of the domain in the $xz$-plane.}\label{fig:fields-2obs-re240}
\end{minipage}

\FloatBarrier%
\section{Error Plots of the 3D Setting}
In Section~\ref{sec:num-example}, we have presented the numerical results
obtained from the DNN-MG method including plots of the velocity and pressure
error for different test cases (cf.~\ref{fig:vp-errors-1obs}
and~\ref{fig:vp-errors-2obs}), showing the improvement in accuracy. To further
substantiate the advantages of DNN-MG we present spatial representations of the
error distribution over the computational domain. The case with 1 obstacle for a
Reynolds number of $\mathrm{Re}=240$, both
for pressure and velocity, is shown in Figure~\ref{fig:fields-1obs-re240-p-error}
and~\ref{fig:fields-1obs-re240-error}. Analogously, the case with 2 obstacles at
the same Reynolds number is shown in Figure~\ref{fig:fields-2obs-re240-p-error}
and~\ref{fig:fields-2obs-re240-error}.
\subsection*{Channel flow with a single obstacle}
\begin{minipage}{\textwidth}
  \centering
  \begin{tikzpicture}[text=white]
    \node[anchor=south west,inner sep=0] (A) at (0,0)
    {\includegraphics[width=\textwidth]{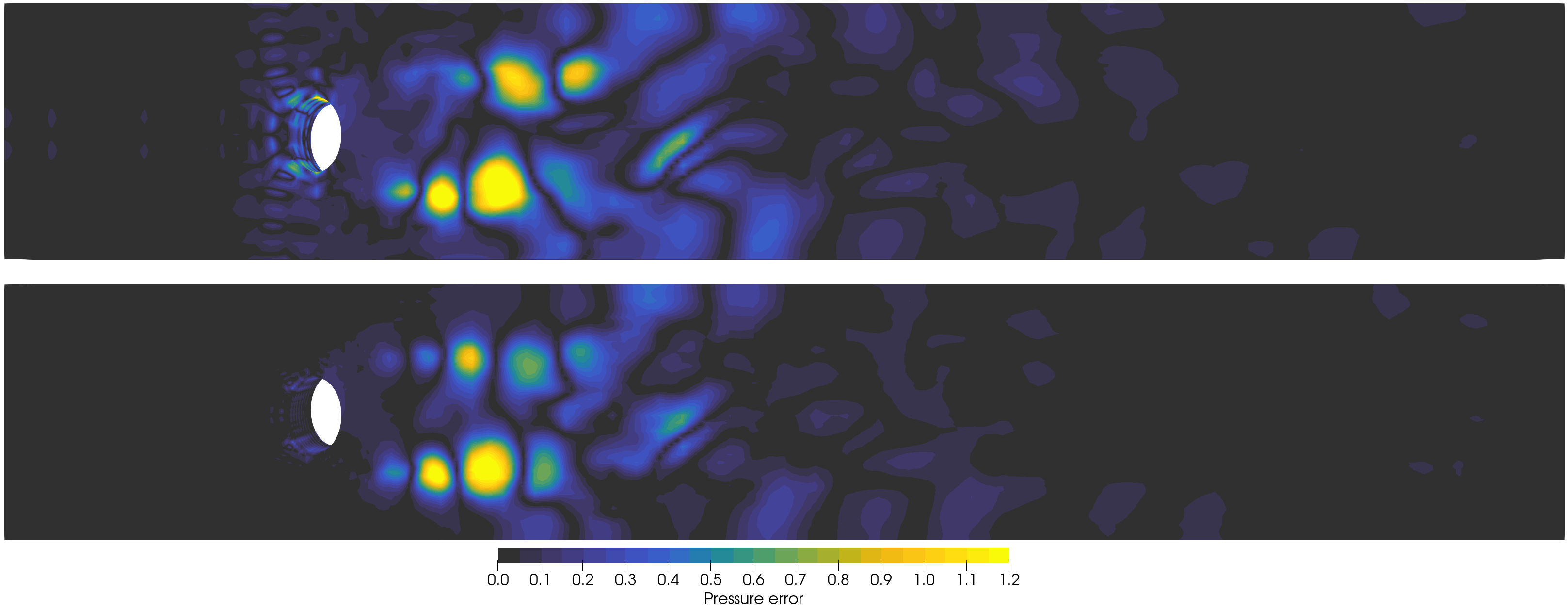}};
    \node[anchor=south west,inner sep=0] (Aa) at (0.3,5.75) {$\text{MG}(3)$};
    \node[anchor=south west,inner sep=0] (Ab) at (0.3,2.9)
    {$\text{DNN-MG}(3+2)$};
    \node[anchor=south west,inner sep=0] (B)
    [below=of A] {
      \includegraphics[width=\textwidth]{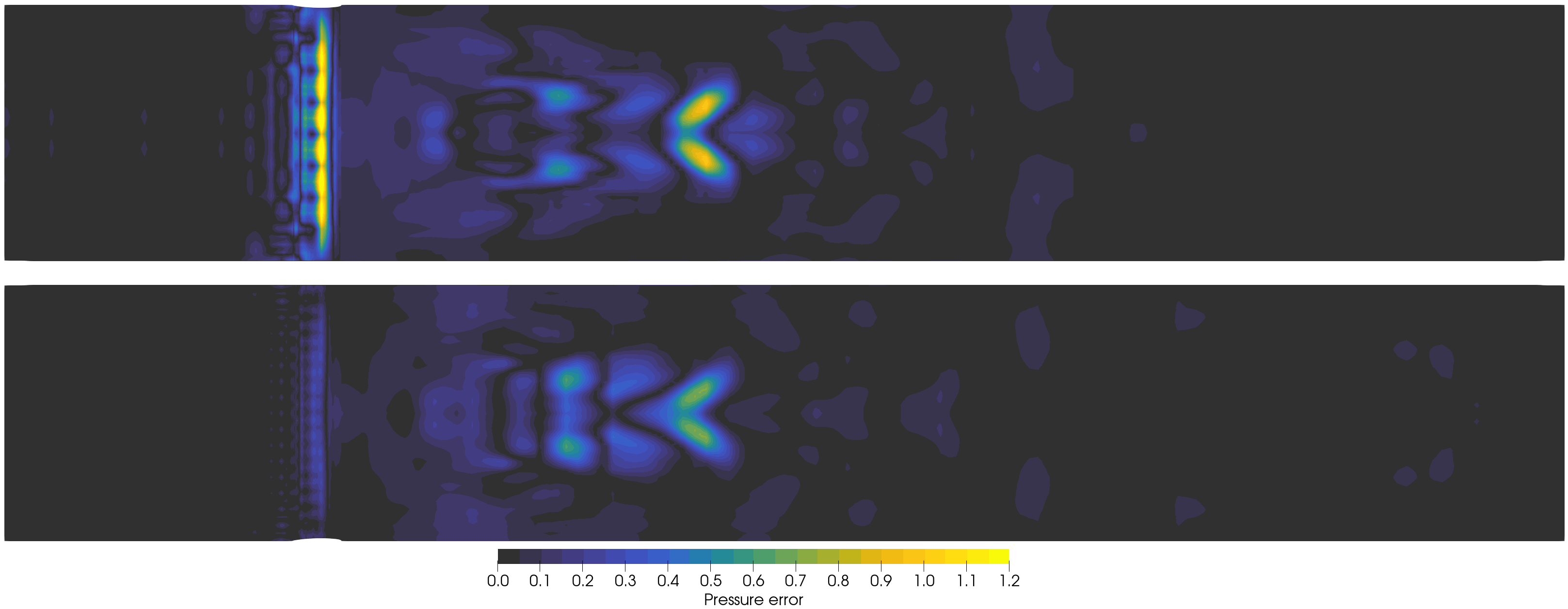}};
    \node[anchor=south west,inner sep=0] (Ba) at (0.3,-1.5) {$\text{MG}(3)$};
    \node[anchor=south west,inner sep=0] (Bb) at (0.3,-4.35)
    {$\text{DNN-MG}(3+2)$};
  \end{tikzpicture}
  \captionof{figure}{Comparison of the error of the pressure for the channel
    with one cylindrical obstacle at time $t=7.5$ at Reynolds number
    $\mathrm{Re}=240$ (cf.\ Figure~\ref{fig:fields-1obs-re240}). The
    results are shown for $\text{DNN-MG}(3+2)$, as well as the coarse
    mesh solution at level $3$ ($\text{MG}(3)$) with respect to the reference with $5$ levels ($\text{MG}(5)$). The
    first two images show a cross-section in the $xy$-plane and
    the others show a cross-section in the $xz$-plane.}\label{fig:fields-1obs-re240-p-error}
\end{minipage}\\

We have seen great improvement in terms of the drag and lift coefficients in
Section~\ref{sec:num-example}. In accordance with this, we see that the
error around the obstacles is greatly reduced by DNN-MG.\@ However, we observe a
significant error reduction over the entire domain. The main error contributions
are due to the different flow patterns evolving in the wake of the channel,
which we can also observe in Figures~\ref{fig:fields-1obs-re240}
and~\ref{fig:fields-2obs-re240}. Towards the end of the channels these
differences seem to dissipate. Overall the spatial distribution of the error
shows that DNN-MG is capable of globally improving the errors via the local
approach. Through the localized approach the method could be extended to use
localized error estimators to train the network on the critical regions in the
combination with an online learning approach.\\

\begin{minipage}{\textwidth}
  \centering
  \begin{tikzpicture}[text=white]
    \node[anchor=south west,inner sep=0] (A) at (0,0)
    {\includegraphics[width=\textwidth]{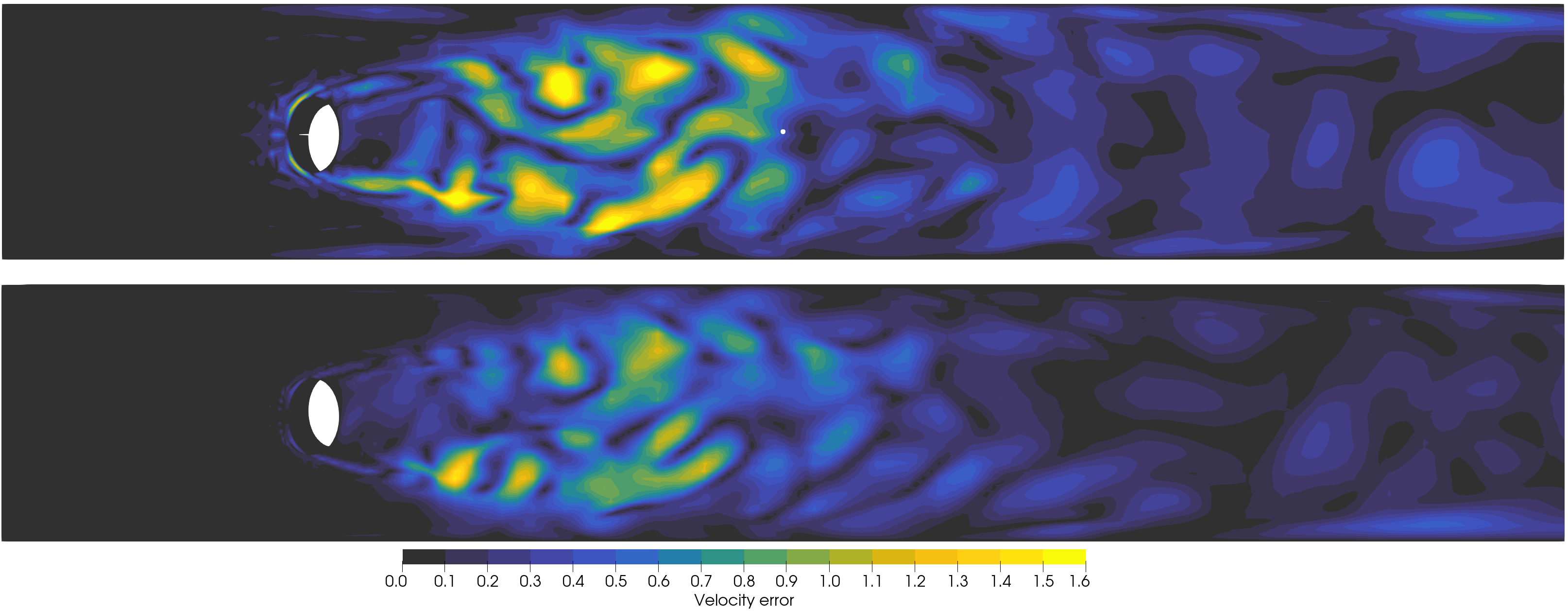}};
    \node[anchor=south west,inner sep=0] (Aa) at (0.3,5.75) {$\text{MG}(3)$};
    \node[anchor=south west,inner sep=0] (Ab) at (0.3,2.9)
    {$\text{DNN-MG}(3+2)$};
    \node[anchor=south west,inner sep=0] (B)
    [below=of A] {
      \includegraphics[width=\textwidth]{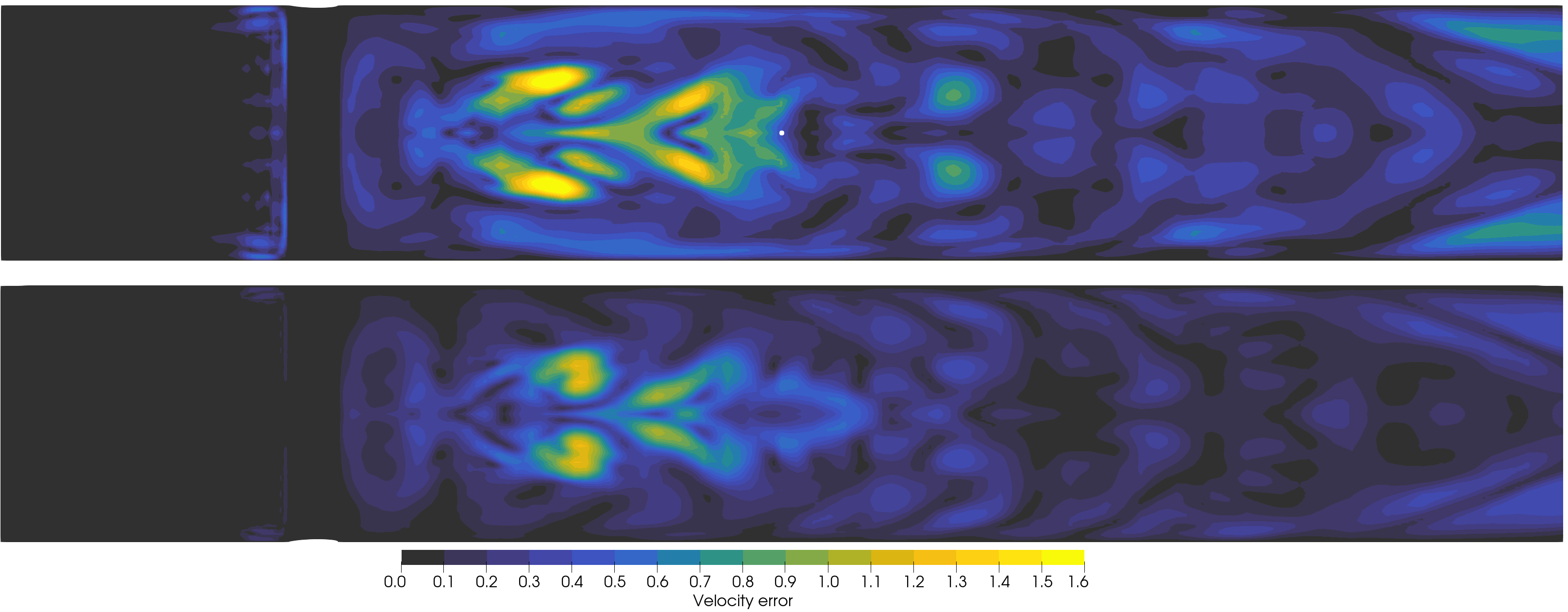}};
    \node[anchor=south west,inner sep=0] (Ba) at (0.3,-1.5) {$\text{MG}(3)$};
    \node[anchor=south west,inner sep=0] (Bb) at (0.3,-4.35)
    {$\text{DNN-MG}(3+2)$};;
  \end{tikzpicture}
  \captionof{figure}{Comparison of the error of the velocity field for the channel
    with one cylindrical obstacle at time $t=7.5$ at Reynolds number
    $\mathrm{Re}=240$ (cf.\ Figure~\ref{fig:fields-1obs-re240}). The
    results are shown for $\text{DNN-MG}(3+2)$, as well as the coarse
    mesh solution at level $3$ ($\text{MG}(3)$) with respect to the reference with $5$ levels ($\text{MG}(5)$). The
    first two images show a cross-section in the $xy$-plane and
    the others show a cross-section in the $xz$-plane.}\label{fig:fields-1obs-re240-error}
\end{minipage}

\subsection*{Channel flow with two obstacles}
\begin{minipage}{\textwidth}
  \centering
  \begin{tikzpicture}[text=white]
    \node[anchor=south west,inner sep=0] (A) at (0,0)
    {\includegraphics[width=\textwidth]{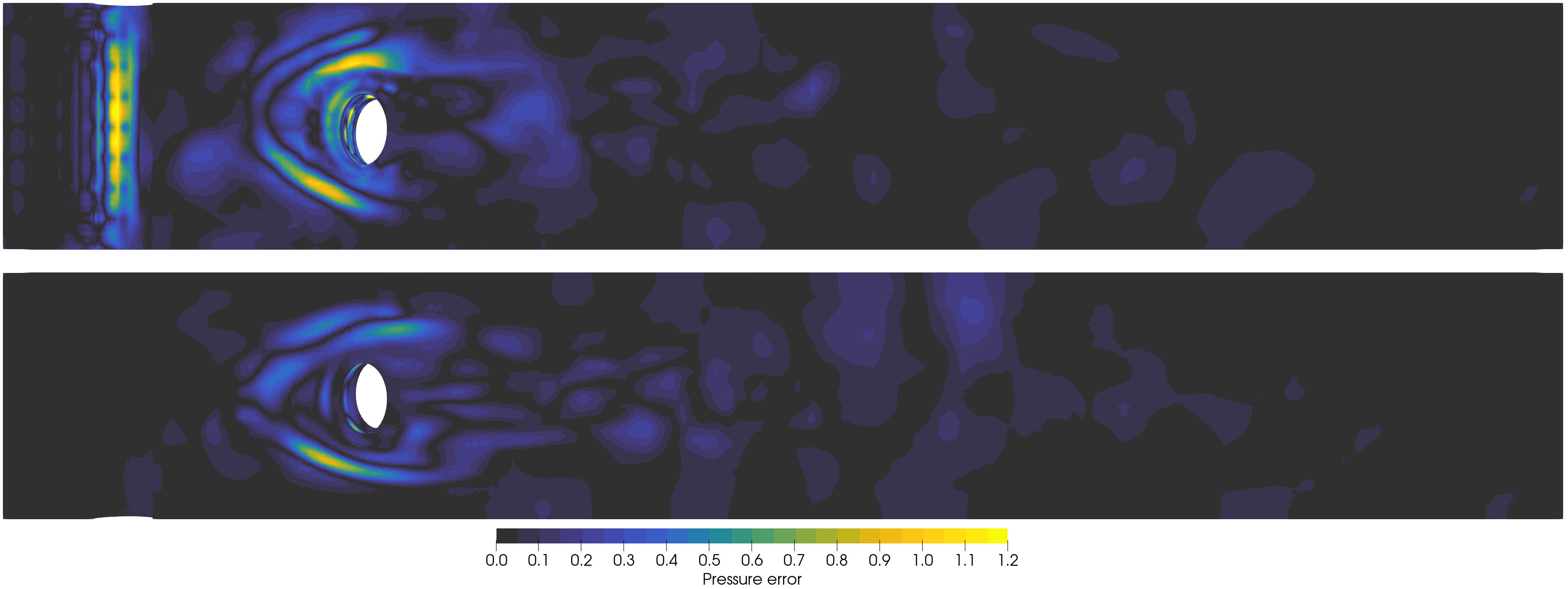}};
    \node[anchor=south west,inner sep=0] (Aa) at (0.3,5.6) {$\text{MG}(3)$};
    \node[anchor=south west,inner sep=0] (Ab) at (0.3,2.8)
    {$\text{DNN-MG}(3+2)$};
    \node[anchor=south west,inner sep=0] (B)
    [below=of A] {
      \includegraphics[width=\textwidth]{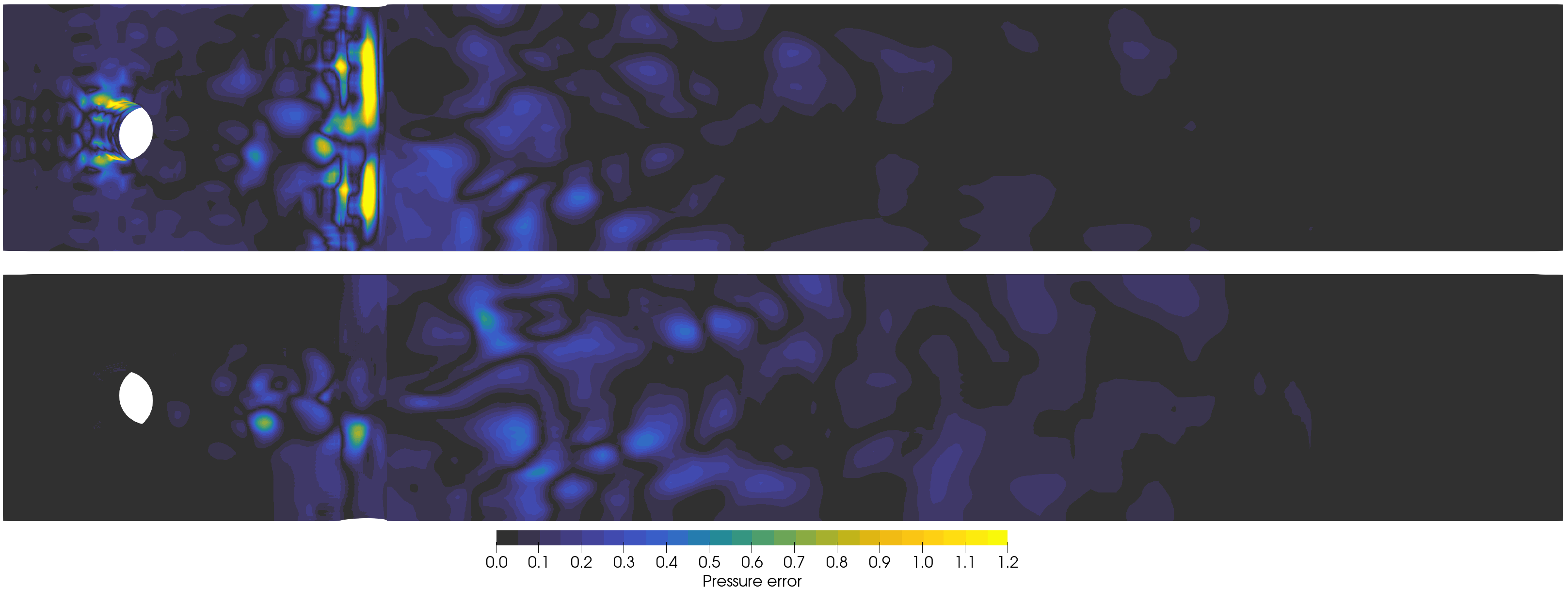}};
    \node[anchor=south west,inner sep=0] (Ba) at (0.3,-1.55) {$\text{MG}(3)$};
    \node[anchor=south west,inner sep=0] (Bb) at (0.3,-4.25)
    {$\text{DNN-MG}(3+2)$};;
  \end{tikzpicture}
  \captionof{figure}{Comparison of the error of the pressure for the channel
    with two cylindrical obstacles at time $t=7.5$ at Reynolds number
    $\mathrm{Re}=240$ (cf.\ Figure~\ref{fig:fields-2obs-re240}). The
    results are shown for $\text{DNN-MG}(3+2)$, as well as the coarse
    mesh solution at level $3$ ($\text{MG}(3)$) with respect to the reference with $5$ levels ($\text{MG}(5)$). The
    first two images show a cross-section in the $xy$-plane and
    the others show a cross-section in the $xz$-plane.}\label{fig:fields-2obs-re240-p-error}
\end{minipage}

\begin{minipage}{\textwidth}
  \centering
  \begin{tikzpicture}[text=white]
    \node[anchor=south west,inner sep=0] (A) at (0,0)
    {\includegraphics[width=\textwidth]{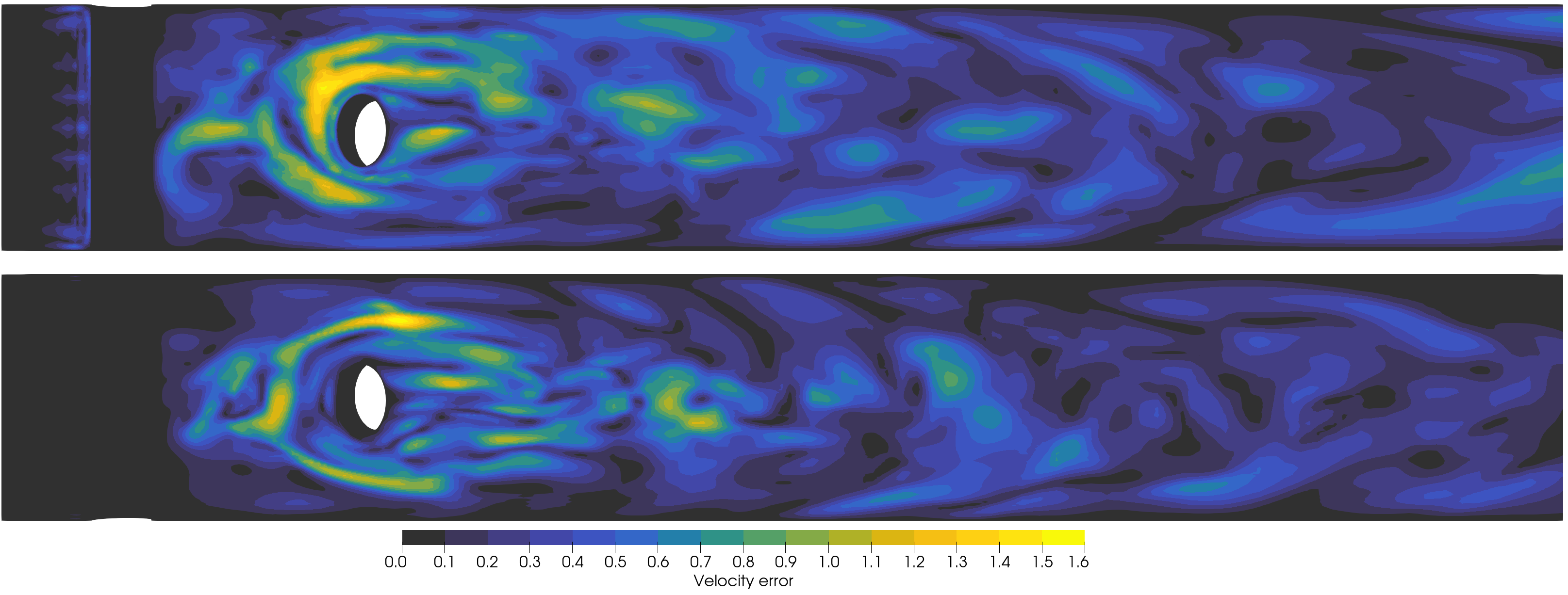}};
    \node[anchor=south west,inner sep=0] (Aa) at (0.3,5.6) {$\text{MG}(3)$};
    \node[anchor=south west,inner sep=0] (Ab) at (0.3,2.8)
    {$\text{DNN-MG}(3+2)$};
    \node[anchor=south west,inner sep=0] (B)
    [below=of A] {
      \includegraphics[width=\textwidth]{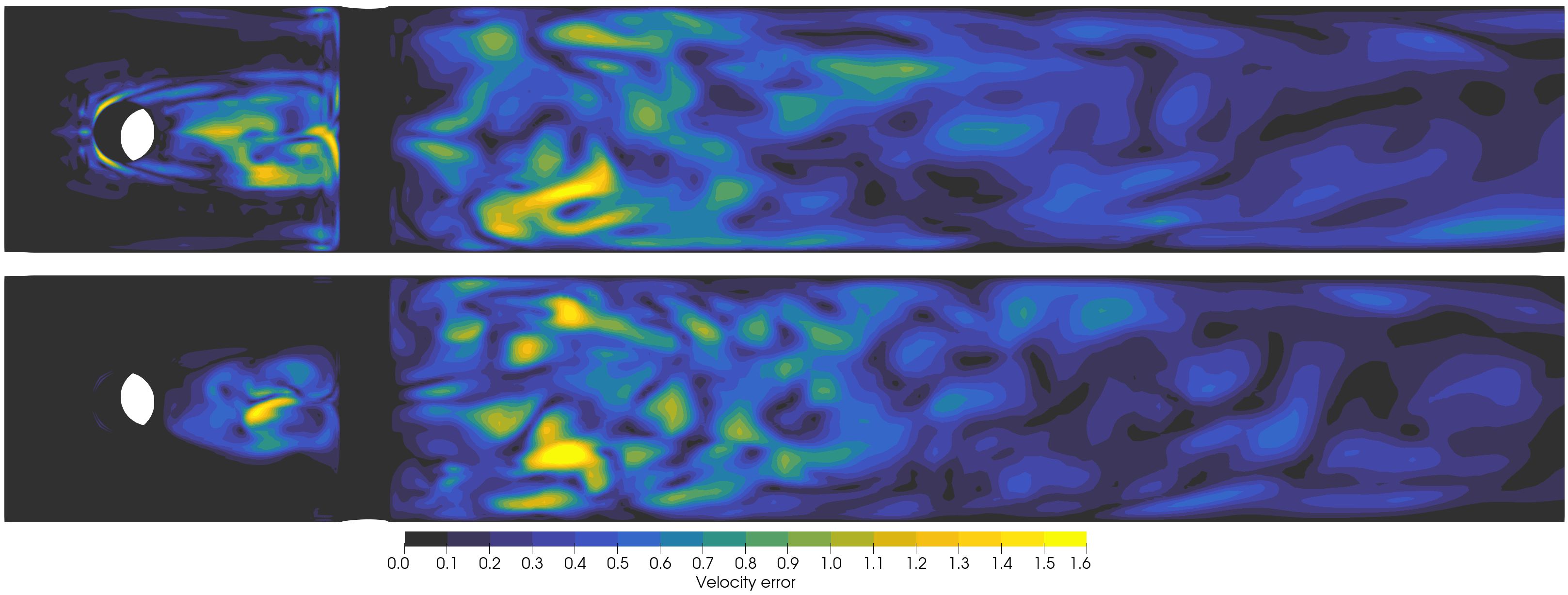}};
    \node[anchor=south west,inner sep=0] (Ba) at (0.3,-1.55) {$\text{MG}(3)$};
    \node[anchor=south west,inner sep=0] (Bb) at (0.3,-4.25)
    {$\text{DNN-MG}(3+2)$};;
  \end{tikzpicture}
  \captionof{figure}{Comparison of the error of the velocity field for the channel
    with two cylindrical obstacles at time $t=7.5$ at Reynolds number
    $\mathrm{Re}=240$ (cf.\ Figure~\ref{fig:fields-2obs-re240}). The
    results are shown for $\text{DNN-MG}(3+2)$, as well as the coarse
    mesh solution at level $3$ ($\text{MG}(3)$) with respect to the reference with $5$ levels ($\text{MG}(5)$). The
    first two images show a cross-section in the $xy$-plane and
    the others show a cross-section in the $xz$-plane.}\label{fig:fields-2obs-re240-error}
\end{minipage}
\end{document}